\documentclass[11pt]{amsart}

\usepackage{xcolor}
\usepackage[all,color]{xy}
\usepackage{amsmath,amsthm, amscd, amssymb, amsfonts}

\usepackage{tikz} 
\usepackage{tikz-cd}

\newcommand{\coev}{\mbox{coev}}
\newcommand{\ev}{\mbox{ev}}

\newcommand{\otk}{{\otimes_{\ku}}}

\newcommand{\Mo}{{\mathcal M}}
\newcommand{\No}{{\mathcal N}}

\newcommand{\uc}{f}

\newcommand{\ca}{{\mathcal C}}

\newcommand{\ot}{{\otimes}}

\newcommand{\op}{\rm{op}}

\newcommand{\Ac}{{\mathcal A}}
\newcommand{\Sc}{{\mathcal S}}

\newcommand\Comod{\operatorname{Comod}}

\newcommand{\Do}{{\mathcal D}}

\newcommand{\vc}{{\mathcal C}}
\newcommand{\Bc}{{\mathcal B}}
\newcommand{\B}{{\mathcal B}}
\newcommand{\ac}{{\mathcal A}}

\newcommand{\YD}{{\mathcal YD}}

\newcommand{\rev}{\rm{rev}}

\newcommand{\ra}{\rm{ra}}
\newcommand{\la}{\rm{la}}

\newcommand{\ku}{{\Bbbk}}

\newcommand{\uno}{ \mathbf{1}}
\newcommand{\C}{{\mathcal C}}

\newcommand{\id}{\mbox{\rm id\,}}

\newcommand{\Id}{\mbox{\rm Id\,}}

\newcommand{\Vect}{\mbox{\rm Vect\,}}
\newcommand{\vect}{\mbox{\rm vect\,}}

\newcommand{\Nat}{\mbox{\rm Nat\,}}
\newcommand{\Rex}{\mbox{\rm Rex\,}}
\newcommand{\Fun}{\operatorname{Fun}}
\newcommand\Rep{\operatorname{Rep}}

\makeatletter
\@namedef{subjclassname@2020}{\textup{2020} Mathematics Subject Classification}
\makeatother

\newcommand\Hom{\operatorname{Hom}}
\newcommand\uhom{\underline{\Hom}}

\renewcommand{\_}[1]{\mbox{$_{\left( #1 \right)}$}}

\theoremstyle{plain}

\numberwithin{equation}{section}

\newtheorem{question}{Question}[section]

\newtheorem{teo}{Theorem}[section]

\newtheorem{lema}[teo]{Lemma}

\newtheorem{cor}[teo]{Corollary}

\newtheorem{prop}[teo]{Proposition}

\newtheorem{claim}{Claim}[section]

\theoremstyle{definition}

\newtheorem{defi}[teo]{Definition}

  \newtheorem{exa}[teo]{Example}

\theoremstyle{remark}

\newtheorem{rmk}[teo]{Remark}

\def\pf{\begin{proof}}

\def\epf{\end{proof}}

\usepackage{comment}

    \newcommand{\pentagon}[7]{ 
    \def\UpSpaceA{$#1$}%
    \def\UpMorAB{$#2$}%
    \def\UpSpaceB{$#3$}%
    \def\UpMorBC{$#4$}%
    \def\UpSpaceC{$#5$}%
    \def\UpMorCD{$#6$}%
    \def\UpSpaceD{$#7$}%
    \pentagoncontinued
    }
    \newcommand{\pentagoncontinued}[3]{ 
    \def\DownMorAE{$#1$}%
    \def\DownSpaceE{$#2$}%
    \def\DownMorED{$#3$}%
          \begin{tikzpicture}[commutative diagrams/every diagram]
          \node (P1) at (180:5cm)
          {\UpSpaceA};
          \node (P2) at (180-45:4cm)
          {\UpSpaceB};
          \node (P3) at (180-45-90:4cm)
          {\UpSpaceC} ;
          \node (P4) at (180-45-90-45:5cm)
          {\UpSpaceD};
          \node (P0) at (180+90:2cm)
          {\DownSpaceE};
          \path[commutative diagrams/.cd, every arrow, every label]
          (P1) edge node[swap] {\UpMorAB} (P2)
          (P2) edge node[swap] {\UpMorBC} (P3)
          (P3) edge node[swap] {\UpMorCD} (P4)
          (P1) edge node {\DownMorAE} (P0)
          (P0) edge node {\DownMorED} (P4);
          \end{tikzpicture}
    } 

    \newcommand{\pentagonABxBCxCDxDExAE}[7]{ 
    \def\SpaceA{$#1$}%
    \def\MorAB{$#2$}%
    \def\SpaceB{$#3$}%
    \def\MorBC{$#4$}%
    \def\SpaceC{$#5$}%
    \def\MorCD{$#6$}%
    \def\SpaceD{$#7$}%
    \pentagonABxBCxCDxDExAEcontinued
    }
    \newcommand{\pentagonABxBCxCDxDExAEcontinued}[3]{ 
    \def\MorDE{$#1$}%
    \def\SpaceE{$#2$}%
    \def\MorEA{$#3$}%
          \begin{tikzpicture}[commutative diagrams/every diagram]
          \node (A) at (180:5cm)
          {\SpaceA};
          \node (B) at (180-45:4cm)
          {\SpaceB};
          \node (C) at (180-45-90:4cm)
          {\SpaceC} ;
          \node (D) at (180-45-90-45:5cm)
          {\SpaceD};
          \node (E) at (180+90:2cm)
          {\SpaceE};
          \path[commutative diagrams/.cd, every arrow, every label]
          (A) edge node[swap] {\MorAB} (B)
          (B) edge node[swap] {\MorBC} (C)
          (C) edge node[swap] {\MorCD} (D)
          (D) edge node[swap] {\MorDE} (E)
          (A) edge node {\MorEA} (E);
          \end{tikzpicture}
}

    \newcommand{\pentagonBAxBCxCDxDExAE}[7]{ 
    \def\SpaceA{$#1$}%
    \def\MorBA{$#2$}%
    \def\SpaceB{$#3$}%
    \def\MorBC{$#4$}%
    \def\SpaceC{$#5$}%
    \def\MorCD{$#6$}%
    \def\SpaceD{$#7$}%
    \pentagonBAxBCxCDxDExAEcontinued
    }
    \newcommand{\pentagonBAxBCxCDxDExAEcontinued}[3]{ 
    \def\MorDE{$#1$}%
    \def\SpaceE{$#2$}%
    \def\MorEA{$#3$}%
          \begin{tikzpicture}[commutative diagrams/every diagram]
          \node (A) at (180:5cm)
          {\SpaceA};
          \node (B) at (180-45:4cm)
          {\SpaceB};
          \node (C) at (180-45-90:4cm)
          {\SpaceC} ;
          \node (D) at (180-45-90-45:5cm)
          {\SpaceD};
          \node (E) at (180+90:2cm)
          {\SpaceE};
          \path[commutative diagrams/.cd, every arrow, every label]
          (B) edge node {\MorBA} (A)
          (B) edge node[swap] {\MorBC} (C)
          (C) edge node[swap] {\MorCD} (D)
          (D) edge node[swap] {\MorDE} (E)
          (A) edge node {\MorEA} (E);
          \end{tikzpicture}
}

\newcommand{\pentagonC}[7]{ 
    \def\SpaceA{$#1$}%
    \def\MorAB{$#2$}%
    \def\SpaceB{$#3$}%
    \def\MorBC{$#4$}%
    \def\SpaceC{$#5$}%
    \def\MorCD{$#6$}%
    \def\SpaceD{$#7$}%
    \pentagonCcontinued
    }
    \newcommand{\pentagonCcontinued}[3]{ 
    \def\MorDE{$#1$}%
    \def\SpaceE{$#2$}%
    \def\MorEA{$#3$}%
          \begin{tikzpicture}[commutative diagrams/every diagram]
          \node (A) at (180:5cm)
          {\SpaceA};
          \node (B) at (180-45:4cm)
          {\SpaceB};
          \node (C) at (180-45-90:4cm)
          {\SpaceC} ;
          \node (D) at (180-45-90-45:5cm)
          {\SpaceD};
          \node (E) at (180+90:2cm)
          {\SpaceE};
          \path[commutative diagrams/.cd, every arrow, every label]
          (A) edge node[swap] {\MorAB} (B)
          (C) edge node {\MorBC} (B)
          (D) edge node {\MorCD} (C)
          (E) edge node {\MorDE} (D)
          (E) edge node[swap] {\MorEA} (A);
          \end{tikzpicture}
}

\newcommand{\pentagonD}[7]{ 
    \def\SpaceA{$#1$}%
    \def\MorBA{$#2$}%
    \def\SpaceB{$#3$}%
    \def\MorBC{$#4$}%
    \def\SpaceC{$#5$}%
    \def\MorCD{$#6$}%
    \def\SpaceD{$#7$}%
    \pentagonDcontinued
    }
    \newcommand{\pentagonDcontinued}[3]{ 
    \def\MorDE{$#1$}%
    \def\SpaceE{$#2$}%
    \def\MorEA{$#3$}%
          \begin{tikzpicture}[commutative diagrams/every diagram]
          \node (A) at (180:5cm)
          {\SpaceA};
          \node (B) at (180-45:4cm)
          {\SpaceB};
          \node (C) at (180-45-90:4cm)
          {\SpaceC} ;
          \node (D) at (180-45-90-45:5cm)
          {\SpaceD};
          \node (E) at (180+90:2cm)
          {\SpaceE};
          \path[commutative diagrams/.cd, every arrow, every label]
          (B) edge node {\MorBA} (A)
          (C) edge node {\MorBC} (B)
          (D) edge node {\MorCD} (C)
          (E) edge node {\MorDE} (D)
          (E) edge node[swap] {\MorEA} (A);
          \end{tikzpicture}
}

\usepackage{calc}
\newcommand{\grau}{gray!60}
\newenvironment{grform}{\begin{tikzpicture}[intext]}{\end{tikzpicture}}
\tikzset{
norm/.style = {ultra thick, color = #1},
norm/.default = black,
intext/.style = {baseline = (current bounding box.center)},
}

\newcommand{\vLine}[5]{
  \draw[ultra thick, color = #5, rounded corners] (#1 , #2) -- (#1 , #2 * 0.7 + #4 * 0.3 ) -- ( #3 , #2 * 0.3 + #4 * 0.7 ) -- (#3 ,#4);
}
\newcommand{\vLineO}[5]{
  \draw[line width = 6pt , color = white, rounded corners] (#1 , #2) -- (#1 , #2 * 0.7 + #4 * 0.3 ) -- ( #3 , #2 * 0.3 + #4 * 0.7 ) -- (#3 ,#4);
  \draw[ultra thick, color = #5, rounded corners] (#1 , #2) -- (#1 , #2 * 0.7 + #4 * 0.3 ) -- ( #3 , #2 * 0.3 + #4 * 0.7 ) -- (#3 ,#4);
}

\newcommand{\dMult}[5]{
\draw[ultra thick, color = #5] (#1 , #2) -- (#1 , #2 + #4 * 0.2) .. controls (#1
, #2 + #4 * 0.5) .. (#1 + #3 *0.25 , #2 + #4 * 0.5) 
-- (#1 + #3 *0.5 , #2 + #4 * 0.5) -- (#1 + #3 *0.75 , #2 + #4 * 0.5) .. controls
(#1 + #3 , #2 + #4 * 0.5) .. (#1 + #3 , #2 + #4 * 0.2) -- (#1 + #3 , #2);
\draw[ultra thick, color = #5] (#1 + #3 *0.5 , #2 + #4 * 0.5) -- (#1 + #3 *0.5 ,
#2 + #4);
}

\newcommand{\dMultO}[5]{
\draw[line width = 6pt , color = white] (#1 , #2) -- (#1 , #2 + #4 * 0.2) ..
controls (#1 , #2 + #4 * 0.5) .. (#1 + #3 *0.25 , #2 + #4 * 0.5) 
-- (#1 + #3 *0.5 , #2 + #4 * 0.5) -- (#1 + #3 *0.75 , #2 + #4 * 0.5) .. controls
(#1 + #3 , #2 + #4 * 0.5) .. (#1 + #3 , #2 + #4 * 0.2) -- (#1 + #3 , #2);
\draw[ultra thick, color = #5] (#1 , #2) -- (#1 , #2 + #4 * 0.2) .. controls (#1
, #2 + #4 * 0.5) .. (#1 + #3 *0.25 , #2 + #4 * 0.5) 
-- (#1 + #3 *0.5 , #2 + #4 * 0.5) -- (#1 + #3 *0.75 , #2 + #4 * 0.5) .. controls
(#1 + #3 , #2 + #4 * 0.5) .. (#1 + #3 , #2 + #4 * 0.2) -- (#1 + #3 , #2);
\draw[ultra thick, color = #5] (#1 + #3 *0.5 , #2 + #4 * 0.5) -- (#1 + #3 *0.5 ,
#2 + #4);
}

\newcommand{\dAction}[6]{
\draw[ultra thick, color = #5] (#1 , #2) -- (#1 , #2 + #4 * 0.2) .. controls (#1
, #2 + #4 * 0.5) .. (#1 + #3 * 0.8 , #2 + #4 * 0.5) -- (#1 + #3 , #2 + #4 *
0.5);
\draw[ultra thick, color = #6] (#1 + #3 , #2) -- (#1 + #3 , #2 + #4);
}

\usepackage{hyperref}

\definecolor{darkgreen}{rgb}{0,0.5,0}
\definecolor{darkblue}{rgb}{0,0.1,0.5}

\usepackage{tikz-cd}

\hypersetup{
    colorlinks=true, 
    linkcolor=darkblue, 
    urlcolor=blue, 
    linktoc=all, 
    citecolor=darkgreen
    }

\newtheoremstyle{introTheorems}
  {\topsep}
  {\topsep}
  {\itshape}
  {0pt}
  {\bfseries}
  {}
  { }
  {\thmname{#1}
  \textnormal{\thmnote{#3}.}
  }

\theoremstyle{introTheorems}

\title[Fibre functors and reconstruction of Hopf algebras]
{Fibre functors and reconstruction of Hopf algebras}
\author[    Lentner and Mombelli]{ Simon Lentner and Mart\'in Mombelli
 }

\begin{document}
\address{Universität Hamburg
\newline \indent
Bundesstrassee 55, D – 20 
\newline \indent
146 Hamburg, Deutschland}

\email{simon.lentner@uni-hamburg.de
\newline \indent\emph{URL:}\/ 
https://www.math.uni-hamburg.de/home/lentner
}

\address{Facultad de Matem\'atica, Astronom\'\i a y F\'\i sica
\newline \indent
Universidad Nacional de C\'ordoba
\newline
\indent CIEM -- CONICET
\newline \indent Medina Allende s/n
\newline
\indent (5000) Ciudad Universitaria, C\'ordoba, Argentina}
\email{martin10090@gmail.com, martin.mombelli@unc.edu.ar
\newline \indent\emph{URL:}\/ https://www.famaf.unc.edu.ar/$\sim$mombelli}

\begin{abstract}
The main objective of the present paper is to present a version of the Tannaka-Krein type reconstruction  Theorems:  If $F:\Bc\to \vc$ is an exact faithful monoidal functor of tensor categories, one would like to realize $\Bc$ as  category of representations of a braided Hopf algebra $H(F)$ in $\vc$. We prove that this is the case iff $\Bc$ has the additional structure of a monoidal $\vc$-module category compatible with $F$, which equivalently means that $F$ admits a monoidal section. For Hopf algebras, this reduces to a version of the Radford projection theorem. The Hopf algebra is constructed through the  relative coend for module categories. We expect this basic result to have a wide range of applications, in particular in the absence of fibre functors, and we give some applications. One particular motivation was the logarithmic Kazhdan-Lusztig conjecture. 
\end{abstract}

 \date{\today}
\maketitle

\tableofcontents

\section*{Introduction}

\subsection*{Background}

The main objective of the present paper is to present a version of the Tannaka-Krein type reconstruction theorems. It is known that if $\Bc$ is a $\ku$-linear abelian locally finite category, $F:\Bc\to \vect_\ku$ is a fibre functor, that is an exact faithful functor, then the coend
$$ C=\int^{B\in \Bc} F(B)\otk F(B)^*,$$
has structure of $\ku$-coalgebra, and the functor $F$ factorizes as 
\begin{equation*}
\xymatrix{
\Bc
\ar[rrr]^-{\qquad\widehat{F}}
\ar[drr]_-{F}
&&&
\mathrm{Comod}(C)
\ar[dl]^-{\uc}\ar@{}[d]|{}&\\
&&\vect_\ku,&
}
\end{equation*} 
where $f:\mathrm{Comod}(C)\to \vect_\ku$ is the forgetful functor, and $\widehat{F}$ is a category equivalence. We refer the reader to \cite{H}, \cite{PS} and references therein. 
If, moreover, $\Bc$ is a monoidal rigid category, and $F$ is a monoidal functor, then it is possible to endow the coalgebra $C$ with structure of Hopf algebra. Succinctly, monoidality of $\Bc$ gives $C$ a product, turning it into a bialgebra, and rigidity of $\Bc$ endows $C$ with an antipode. 

\medbreak

Some generalizations of this result appeared in the literature. We only mention some of them. In these more general versions, a fibre functor $F:\Bc\to \vc$ is considered, where $\vc$ is an arbitrary monoidal category. In some versions, the object $C$ is not a Hopf algebra in $\ca$. For example in \cite{BV},  it is shown that if the functor $F$ has a right adjoint $G:\ca\to \Bc$, then, the associated monad to this adjunction $T=F\circ G:\vc\to \vc$ is a bimonad, and there is a commutative diagram
\begin{equation*}
\vspace{.1cm}
\xymatrix{
\Bc\;\;
\ar[rr]^-{\widehat{F}}
\ar[dr]_-{F}
&&
\mathcal{C}^T
\ar[dl]^-{\uc}\ar@{}[d]|{}&\\
&\;\mathcal{C}.
\ar@/^1.5pc/@[gray][ul]^{G\;\;}&
}
\end{equation*} 

Here $f: \vc^T\to \vc$ is the forgeful functor. Beck's monadicity Theorem implies that the functor $\widehat{F}:\Bc\to \ca^T$ is a monoidal equivalence. Rigidity of both monoidal categories $\Bc, \vc$ imply that $T$ has a Hopf monad structure. In \cite{Li1}, \cite{Li2} Lyubashenko reconstructed the object $C$ as a coend $\int^{B\in \Bc} F(B)\boxtimes F(B)^*$, belonging to some completion of the Deligne tensor product $\vc \boxtimes \vc$, and it turns out to be a {\it squared coalgebra}. In the work of Majid \cite{Majid}, he started with a monoidal functor $F:\Bc\to \vc$, where $\vc$ is a braided monoidal category, and he reconstructed a Hopf algebra $C=\int^{B\in \Bc} F(B)^*\ot F(B)$, and set up a commutative diagram
\begin{equation*}
\vspace{.1cm}
\xymatrix{
\Bc\;\;
\ar[rr]^-{\widehat{F}}
\ar[dr]_-{F}
&&
{}^C\mathcal{C}
\ar[dl]^-{\uc}\ar@{}[d]|{}&\\
&\;\mathcal{C}.&
}
\end{equation*}
With this generality, the functor $\widehat{F}$ no longer need to be an equivalence.


\subsection*{Our approach}

 The point of view of this paper owes a lot to \cite{PS}. However, we shall work only with finite categories. Let $\ku$ be an arbitrary field, and let $\Bc, \vc$ be finite $\ku$-linear abelian categories. Assume both $\ca, \Bc$ are  rigid monoidal categories, and $F:\Bc\to \vc$ is a monoidal functor. If $F$ has a {\it section} $G:\vc\to \Bc$, that is $F\circ G\simeq \Id_\vc$ as monoidal functors, then it is possible to endow $\Bc$ with an action of $\vc$, such that $F$ is a monoidal $\vc$-module functor. This action behaves well together with the monoidal product of $\Bc$, in a sense that we call \textit{monoidal module category}.

  Under these conditions, we aim at constructing a Hopf algebra $H\in \vc$ and obtaining a kind of Radford projection Theorem in this categorical setting.

 \medbreak
 
The reconstruction of such Hopf algebra is given in some steps, that we describe as follows. If $\Bc$ is a right $\vc$-module category, $F:\Bc\to \vc$ is an exact faithful  module functor, we construct a coalgebra $C(F)\in \vc$ as
 $$ C(F)=\oint^{B\in \Bc} F(B)\ot {}^*F(B).$$
 Here $\oint$ stands for the \emph{relative coend}, a new tool, developed in \cite{BM}, in the context of module categories. This tool is one of the new features that we incorporate in these reconstruction theorems. The coproduct and counit of  $C(F)$ are defined using universal properties of the dinatural transformations associated with the relative coend. See Proposition \ref{coalg-coend}. Moreover, we show that the functor $F$ factorizes as 
\begin{equation*}
\vspace{.1cm}
\xymatrix{
\Bc\;\;
\ar[rr]^-{\widehat{F}}
\ar[dr]_-{F}
&&
{}^C\vc
\ar[dl]^-{\uc}\ar@{}[d]|{}&\\
&\vc.&
}
\end{equation*} 
Here $f:{}^C\vc\to \vc $ is the forgetful functor, and $\widehat{F}: \Bc\to {}^C\vc$ is an equivalence of $\vc$-module categories. The existence of $\widehat{F}$ is stated in Proposition \ref{factorization-F}, and the proof that it is an equivalence of categories is given in Theorem \ref{recons-th}.

If in addition $\vc$ is a braided rigid monoidal  $\vc$-module category,  $\Bc$ is a rigid monoidal category, and $F$ is a monoidal functor, then we endow the coalgebra $C(F)$ with a product, turning it into a bialgebra in $\vc$. Rigidity of $\Bc$ allows to define an antipode on $C(F)$, making it into a braided Hopf algebra. These results are stated in Theorem \ref{hopf-rec} and Corrollary \ref{antipode-reconstruct}.

\medbreak 

The contents of the paper are the following. In Section \ref{Section:tensor-categories} we give a brief account of some basic facts about tensor categories and module categories that will be used throughout the paper. In Section \ref{Section:mcoends} we review the notion of  \emph{relative (co)end} of a functor  in the setting of module categories over a tensor category $\ca$. This tool was developed by the second author in \cite{BM} as generalization of the usual (co)end. 
Let $\ca$ be a tensor category, $\Mo$ be a left $\vc$-module category and $\Ac$ is some target category. If a functor 
$S:\Mo^{\op}\times \Mo\to \Ac$ has a \textit{prebalancing}, that is natural isomorphisms
$$\beta^X_{M,N}: S(M,X\triangleright N)\to S(X^*\triangleright M,N),$$
then the relative coend 
$$\oint^{M\in \Mo} (S,\beta)$$
is an object in $ \Ac$ equipped with dinatural transformations $\pi_M:S(M,M)\xrightarrow{ . .} \oint^{M\in \Mo} (S,\beta)$, such that it satisfies some extra condition, see for example \eqref{dinat:coen:module:left}, and it is universal with this property. When $\ca=\vect_\ku$, the relative (co)end coincides with the usual (co)end. 

\medbreak

Since relative (co)ends are objects defined by a universal property, they may not exist in general. Section \ref{Section_coendExists} is devoted to prove that all relative coends, used in this work, actually exist. In this Section, it is crucial that all categories, the tensor category $\vc$ and the module category $\Mo$ are finite. We also require that the action $ \vc\times \Mo\to \Mo$ is exact in each variable, allowing us to use \cite[Thm 2.24]{DSS}, that says that there is an equivalence of module categories $\Mo\simeq \vc_A$, for some algebra $A\in \vc.$

\medbreak

In Section \ref{Section_Hopf} we review the definition of a Hopf algebra $H$ in a braided tensor category $\vc$ and its tensor category of comodules ${^H}\vc$. We also review the definition of $H$-$H$-Yetter Drinfeld modules in this setting, which produces a braided tensor category ${^H_H}\YD(\vc)$. Then we introduce a notion of a $\vc$-module category $\Bc$ with a compatible tensor structure in the sense that there is a natural isomorphism 
$$l_{X,B}: X \triangleright B\xrightarrow{\;\simeq\; } (X\triangleright \uno_\Bc)\ot B,$$
for all $X\in\vc,\,B\in\Bc$, satisfying certain axioms. Note that there is a second reasonable notion of a monoidal module category, where the action $X\triangleright-$ is a monoidal functor, and this definition is not equivalent. Under the presence of a monoidal  functor $F:\Bc\to \ca$, compatibility of the monoidal product in $\Bc$ and the action of $\vc$ on $\Bc$, is equivalent to the existence of a section to $F$, that is a monoidal functor $G:\vc\to \Bc$ such that there is a monoidal natural isomorphism $F\circ G\simeq \Id.$

\subsection*{Acknowledgments} The work of both authors was partially done during the \emph{Workshop on Hopf algebras and Tensor categories} held in May 2023, in Marburg, Germany. We thank the organizers, and in particular to Istvan Heckenberger. The attendance of M.M., to this workshop, was possible due to a  grant of the Alexander von Humboldt Foundation, under a Research Group Linkage Programme. Part of the work of M.M. was done during a visit at Hamburg University, M.M.  thanks Christoph Schweigert for the warm hospitality and for providing us with reference \cite{H}. The work of M.M. was partially supported by Secyt-U.N.C., Foncyt and CONICET Argentina. S.L. thanks T. Gannon and T. Creutzig for hospitality at the University of Alberta and the Alexander von Humboldt Foundation for financial support via the Feodor Lynen Fellowship.

\section{Preliminaries}

Throughout this paper, $\ku$ will denote an arbitrary field. 
 We shall denote by $\vect_\ku$ the category of finite dimensional $\ku$-vector spaces. 
 
 A \emph{finite  category}  \cite{EO} is an abelian $\ku$-linear category such that it has only a finite number of isomorphism
classes of simple objects, Hom spaces are finite-dimensional
$\ku$-vector spaces, all objects have finite lenght and every simple object has a projective cover. All these conditions are equivalent to requiring that, the category is equivalent to the category of finite-dimensional representations of a finite-dimensional $\ku$-algebra.

 \medbreak
 
 If $\Mo, \No$ are categories, and $F:\Mo\to \No$ is a functor, we shall denote by $F^{\ra}, F^{\la}:\No\to \Mo$ its right and left adjoint, respectively. We shall denote by $\Mo^{\op}$ the opposite category. If $f:M\to N$ is a morphism in $\Mo$, sometimes we shall denote by $f^{\op}:N\to M$ the same map but understood as a morphism in  $\Mo^{\op}$.

Any abelian $\ku$-linear category $\Mo$ has a canonical action of the category of finite dimensional $\ku$-vector spaces
\begin{equation}\label{vect-action}
 \bullet: \vect_\ku \times \Mo\to \Mo.
\end{equation}
See for example \cite[Lemma 2.2.2]{PS}. Any additive $\ku$-linear functor $F:\Mo\to \No$ between abelian $\ku$-linear categories respects the action of $\vect_\ku$, that is, there are natural isomorphisms
$$d_{V,M}:F(V \bullet M)\to V \bullet F(M),$$
$V\in \vect_\ku$, $M\in \Mo$, satisfying certain axioms.
\medbreak

\medbreak

 From now on, \emph{all categories} will be assumed to be finite  abelian $\ku$-linear categories, and all functors will be additive $\ku$-linear. Here $\ku$ is an arbitrary field. All our proofs work in the presence of an associator, but for simplicity we assume in the presentation that the categories are strict. 

\section{Representations of tensor categories}\label{Section:tensor-categories}

\subsection{Finite tensor categories}

A  finite tensor category $\ca$ is a monoidal rigid category, with simple unit object $\uno\in \ca$. We refer to \cite{EO} for more details on finite tensor categories. Without loss of generality, we shall assume that tensor categories in this work are strict. 

\medbreak

If $\ca, \Do$ are monoidal categories, a monoidal functor is a functor $F:\ca\to \Do$ equipped with natural isomorphisms
$$\xi_{A,B}:F(A)\ot F(B)\to F(A\ot B), $$
such that
\begin{equation}\label{monoidal-functor}
\xi_{A,B\ot C} (\id_{F(A)}\ot \xi_{B,C})=\xi_{A\ot B, C} (\xi_{A,B}\ot \id_{F(C)}). 
\end{equation}
If $(F,\xi), (\widetilde{F}, \widetilde{\xi}):\ca\to \Do$ are two monoidal functors, a monoidal \textit{natural transformation} between $F$ and $\widetilde{F}$ is a natural transformation $\alpha:F\to \widetilde{F}$ such that
\begin{equation}\label{natural-monoidal}
\alpha_{A\ot B} \xi_{A,B}=\widetilde{\xi}_{A,B} (\alpha_A\ot \alpha_B),
\end{equation}
for any $A,B\in \ca$.
\medbreak

If $\ca$ is a category, for any $X\in \ca$ we shall denote by 
$$ \ev_X:X^*\ot X\to \uno, \quad \coev_X:\uno\to X\ot X^*$$
the evaluation and coevaluation. By abuse of notation, we shall also denote by
$$\ev_X:X\ot {}^*X\to \uno,  \quad \coev_X:\uno\to  {}^*X\ot X$$
the evaluation and coevaluation for the left duals. We will use the following basic result.
If $f:X\to Y$ is an isomorphism in $\ca$ then
\begin{align}\label{duality-iso} \begin{split} (\id_{ {}^*X}\ot f )\coev_X=& ({}^*f\ot \id_Y )\coev_Y,\\
\ev_Y (f\ot  \id_Y)=&\ev_X (\id_X\ot {}^*f).
\end{split}
\end{align}
 It is well known that for any pair of objects $X, Y\in \ca$ there are canonical natural isomorphisms
\begin{equation}\label{duals-product-r}
 \phi^r_{X,Y}: Y^* \ot X^* \to(X\ot Y)^*,
\end{equation} 
\begin{equation}\label{duals-product-l}
 \phi^l_{X,Y}: {}^*Y \ot {}^*X \to {}^*(X\ot Y).
\end{equation} 
These isomorphisms allow us to compute coevaluation and evaluation maps of duals; more precisely we shall need the following identity
\begin{equation}\label{coev-of-dual} \ev_{{}^*V}=  {}^*(\coev_X) \phi^l_{{}^*X,X}.
\end{equation} 

We are going to make use, very often, of the canonical natural isomorphisms
\begin{align}\label{duals1}\begin{split} \psi^Z_{X,Y}:\Hom_\Bc(X\ot {}^*Y,Z)\to \Hom_\Bc(X,Z\ot Y),\\
\psi^Z_{X,Y}(f)=(f\ot\id_Y)(\id_X\ot \coev_Y).
\end{split}
\end{align}
And its inverse
\begin{align}\label{duals2}\begin{split} \bar{\psi}^Z_{X,Y}:\Hom_\Bc(X,Z\ot Y)\to \Hom_\Bc(X\ot {}^*Y,Z)\\
\bar{\psi}^Z_{X,Y}(g)=(\id_Z\ot\ev_Y)(g\ot\id_{{}^*Y}).\end{split}
\end{align}

We shall also need the following basic fact.
\begin{lema} For any pair $V,W\in \ca$ the evaluation of $V\ot W$ is given by
\begin{equation}\label{evaluation-tensor-prod} \ev_{V\ot W} = \ev_V(\id_V\ot \ev_W\ot \id_{{}^*V})(\id_{V\ot W} \ot (\phi^l_{V,W})^{-1}).\qed
\end{equation}
\end{lema}

\subsection{Module categories over tensor categories}
 A  left \emph{module} category over  
$\ca$ is a  category $\Mo$ together with a $\ku$-bilinear 
bifunctor $\rhd: \ca \times \Mo \to \Mo$, exact in each variable,  endowed with 
 natural associativity
and unit isomorphisms 
$$m_{X,Y,M}: (X\otimes Y)\triangleright   M \to X\triangleright  
(Y\triangleright M), \ \ \ell_M: \uno \triangleright  M\to M.$$ 
These isomorphisms are subject to the following conditions: 

\begin{equation}\label{left-modulecat1}
\pentagon
{(X \otimes (Y\otimes Z))\triangleright M }
{\rotatebox{60}{$\sim$}}
{((X \otimes Y)\otimes Z)\triangleright M}
{m_{X\otimes Y, Z, M}}
{(X \otimes Y)\triangleright (Z\triangleright M)}
{m_{X, Y, Z\triangleright M}}
{X \triangleright (Y\triangleright (Z\triangleright M))}
{m_{X, Y\otimes Z, M}}
{X \triangleright ((Y\otimes Z)\triangleright M) }
{\mathrm{id}_{X}\triangleright m_{Y,Z, M}}
\end{equation}
for any $X, Y, Z\in\C, M\in\Mo$, as well as 
\begin{equation}\label{left-modulecat2} (\id_{X}\triangleright \ell_M)m_{X,{\bf 1} ,M}=\id.
\end{equation} 
Sometimes we shall also say  that $\Mo$ is a $\ca$-\emph{module category} or a representation of $\ca$.

\medbreak

Let $\Mo$ and $\Mo'$ be a pair of $\C$-module categories. A \emph{ module functor} is a pair $(F,c)$, where  $F:\Mo\to\Mo'$  is a functor equipped with natural isomorphisms
$$c_{X,M}: F(X\triangleright M)\to
X\triangleright F(M),$$ 
for $X\in  \ca$, $M\in \Mo$, such that
\begin{equation}\label{modfunctor1}
\pentagon
{F((X\otimes Y)\triangleright M)}
{F(m_{X,Y,M})}
{F(X\triangleright (Y\triangleright M))}
{c_{X,Y\triangleright M}}
{X\triangleright F(Y\triangleright M)}
{\mathrm{id}_X \triangleright  c_{Y,M}}
{X\triangleright (Y\triangleright F(M))}
{c_{X\otimes Y,M}}
{(X\otimes Y)\triangleright F(M)}
{m_{X,Y,F(M)}}
\end{equation}
for any $X, Y\in\ca$, $M\in \Mo$, as well as
\begin{align}\label{modfunctor2}
\ell_{F(M)} \,c_{\uno ,M} &=F(\ell_{M}).
\end{align}

Module functors are composable, if $\Mo''$ is a $\ca$-module category and
$(G,d): \Mo' \to \Mo''$ is another module functor then the
composition
\begin{equation}\label{modfunctor-comp}
(G\circ F, e): \Mo \to \Mo'', \qquad  e_{X,M} = d_{X,F(M)}\circ
G(c_{X,M}),
\end{equation} is
also a module functor.

\smallbreak  

A \textit{natural module transformation}, between  module functors $(F,c)$ and $(G,d)$, is a 
 natural transformation $\theta: F \to G$ such
that
\begin{gather}
\label{modf-nat} d_{X,M}\theta_{X\triangleright M} =
(\id_{X}\triangleright \theta_{M})c_{X,M},
\end{gather}
 for any $X\in \ca$, $M\in \Mo$. The vector space of natural module transformations will be denoted by $\Nat_{\!m}(F,G)$. Two module functors $F, G$ are \emph{equivalent} if there exists a natural module isomorphism
$\theta:F \to G$. We denote by $\Fun_{\ca}(\Mo, \Mo')$ the category whose
objects are module functors $(F, c)$ from $\Mo$ to $\Mo'$ and arrows module natural transformations. 

\medbreak
Two $\C$-modules $\Mo$ and $\Mo'$ are {\em equivalent} if there exist module functors $F:\Mo\to
\Mo'$, $G:\Mo'\to \Mo$, and natural module isomorphisms
$\Id_{\Mo'} \to F\circ G$, $\Id_{\Mo} \to G\circ F$.

\medbreak

A module category will be called \textit{strict} if isomorphisms $m$ and $l$ are identities. Any module category is equivalent to a strict one. We will often  assume that the module category is strict without further mention.

\medbreak
 A \emph{right module category} over $\ca$
 is a   category $\Mo$ equipped with an exact
bifunctor $\triangleleft:  \Mo\times  \ca\to \Mo$ and natural   isomorphisms 
$$\widetilde{m}_{M, X,Y}: M\triangleleft (X\ot Y)\to (M\triangleleft X) \triangleleft Y, \quad r_M:M\triangleleft \uno\to M$$ such that
\begin{equation}\label{right-modulecat1} 
\pentagon
{M\triangleleft((X\otimes Y)\otimes Z)}
{\rotatebox{60}{$\sim$}}
{M\triangleleft(X\otimes (Y\otimes Z))}
{\widetilde{m}_{M,X ,Y\ot Z }}
{(M\triangleleft X)\triangleleft (Y\otimes Z)}
{\widetilde{m}_{M\triangleleft X, Y ,Z }}
{((M\triangleleft X)\triangleleft Y)\triangleleft Z)}
{\widetilde{m}_{M,X\ot Y ,Z}}
{(M\triangleleft(X\otimes Y))\triangleleft Z}
{\widetilde{m}_{M,X,Y}\triangleleft \mathrm{id}_Z}
\end{equation}
as well as 
\begin{equation}\label{right-modulecat2} (r_M\triangleleft \id_X)  \widetilde{m}_{M,\uno, X}= \id.
\end{equation}

If $\Mo,  \Mo'$ are right $\ca$-modules, a module functor from $\Mo$ to $  \Mo'$ is a pair $(T, d)$ where
$T:\Mo \to \Mo'$ is a  functor and $d_{M,X}:T(M\triangleleft X)\to T(M)\triangleleft X$ are natural  isomorphisms
such that for any $X, Y\in
\ca$, $M\in \Mo$:
\begin{equation}\label{modfunctor11}
\pentagon
{T(M\triangleleft(X\otimes Y))}
{T(\tilde{m}_{M, X, Y})}
{T((M\triangleleft X)\triangleleft Y)}
{d_{M\triangleleft X, Y}}
{T(M\triangleleft X)\triangleleft Y}
{d_{M,X}\triangleleft \id_Y}
{(T(M)\triangleleft X)\triangleleft Y}
{d_{M, X\ot Y}}
{T(M)\triangleleft(X\otimes Y))}
{\tilde{m}_{T(M), X,Y}}
\end{equation}
as well as
\begin{align}\label{modfunctor21}
r_{T(M)} \,d_{ M,\uno} &=T(r_{M}).
\end{align}

\subsection{The internal Hom}\label{subsection:internal-hom} Let $\ca$ be a  tensor category and $\Mo$ be  a left $\C$-module category. For any pair of objects $M, N\in\Mo$, the \emph{internal Hom} is an object $\uhom(M,N)\in \C$ representing the left exact functor $$\Hom_{\Mo}(-\triangleright M,N):\ca^{\op}\to \vect_\ku.$$ This means that, there are natural isomorphisms
$$\Hom_{\ca}(X,\uhom_\Mo(M,N))\simeq \Hom_{\Mo}(X\triangleright M,N).$$
The internal Hom for right $\ca$-module categories is defined similarly.

\medbreak

The next technical result will be needed later. Recall natural isomorphisms $\phi^l$ defined in \eqref{duals-product-l}.

\begin{lema} Let $\Mo$ be a  right $\ca$-module category with action given by $\triangleleft: \Mo\times \ca\to \Mo$. If $(F,c):\Mo\to \ca$ is a module functor, then, for any $B\in \Mo$, $X\in \ca$, the following diagrams commute. 
\begin{equation}\label{evaluation-mod}
\begin{tikzcd}[column sep=-1em]
&
(F(B)\ot X)\ot {^*}(F(B) \ot X) 
\ar{dr}{\mathrm{id}\otimes (\phi^l_{F(B),X})^{-1}}
&
\\
F(B\triangleleft X)\ot {^*}F(B \triangleleft X) 
\ar{ru}{c_{B,X}\ot {^*}c_{B,X}^{-1}}
\arrow[swap]{dr}{\mathrm{ev}_{F(B\triangleleft X )}}
&
&
(F(B)\ot X)\ot ({^*}X\ot {^*}F(B)) 
\ar{dl}{\mathrm{ev}_{F(B)}(\mathrm{id}\ot \mathrm{ev}_X\ot \mathrm{id})}
&
\\
&
\uno
&
\end{tikzcd}
\end{equation}

\begin{equation}\label{coevaluation-mod}
\begin{tikzcd}[column sep=-1em]
&
({^*}X\ot {^*}F(B)) \ot (F(B)\ot X)
\arrow[swap]{dl}{{^*}c_{B,X} \ot c^{-1}_{B,X}}
&
\\
{^*}F(B \triangleleft X)\otimes F(B\triangleleft X) 
&
&
{^*}(X\ot F(B)) \ot (F(B)\ot X)
\arrow[swap]{ul}{\phi^l_{F(B),X}\otimes \mathrm{id}}
&
\\
&
\uno
\arrow{ul}{\mathrm{coev}_{F(B\triangleleft X )}}
\arrow[swap]{ur}{(\mathrm{id}\ot \mathrm{coev}_{F(B)}\ot\mathrm{id}) \mathrm{coev}_X }
&
\end{tikzcd}
\end{equation}
\end{lema} 
\pf The proof follows, by simply checking that morphisms
$$\ev_{F(B\triangleleft X )} = \ev_{F(B)}\big(\id\ot \ev_X\ot \id \big)\big(c_{B,X}\ot ( {}^*c_{B,X}\phi^l_{F(B),X})^{-1}\big), $$
$$\coev_{F(B\triangleleft X )}= \big( {}^*c_{B,X}\phi^l_{F(B),X} \ot c^{-1}_{B,X}\big)\big(\id\ot \coev_{F(B)}\ot\id\big) \coev_X,$$
satisfy rigidity axioms.
\epf

\section{The (co)end for module categories}\label{Section:mcoends}

In this Section we recall the notion of \textit{relative (co)ends}; a tool developed in \cite{BM} in the context of representations of tensor categories, generalizing the well-known notion of (co)ends in category theory.

\medbreak
 Let $\ca$ be a  tensor category and $\Mo$ be a left $\ca$-module category. Assume that $\Ac$ is a category and $S:\Mo^{\op}\times \Mo\to \Ac$ is a functor equipped with natural isomorphisms
\begin{equation} \beta^X_{M,N}: S(M,X\triangleright N)\to S(X^*\triangleright M,N),
\end{equation}
for any $X\in \ca, M,N\in \Mo$. We shall say that $\beta$ is a \textit{prebalancing} of the functor $S$. Sometimes we shall say that it is a $\ca$-\textit{prebalancing}, to emphasize the dependence on $\ca$.

\begin{defi} The \textit{relative end} of the pair $(S,\beta)$ is an object $E\in \Ac$ equipped with dinatural transformations $\pi_M: E\xrightarrow{ . .} S(M,M)$ such that 
\begin{equation}\label{dinat:end:module:left}
\pentagonC
{S(M,M)}
{S(\mathrm{ev}_X^{op}\triangleright\mathrm{id}_M,\mathrm{id}_M)}
{S((X^*\triangleright X)\triangleright M, M)}
{S(m_{X^*,X,M}^{op} ,  \mathrm{id}_M) \vphantom{{{X^X}^X}^X}} 
{S(X^*\triangleright(X\triangleright M), M)}
{\beta^X_{X\triangleright M,M}}
{S(X\triangleright M,X\triangleright M)}
{\pi_{X\triangleright M}}
{E}
{\pi_{ M}}
\end{equation}
for any $X\in \ca, M\in \Mo$, and is universal with this property. This means that, if $\widetilde{E}\in \Ac$ is another object with dinatural transformations $\xi_M:\widetilde{E}\xrightarrow{ . .} S(M,M)$, such that they fulfill \eqref{dinat:end:module:left}, there exists a unique morphism $h:\widetilde{E}\to E$ such that $\xi_M= \pi_M\circ h $.
\end{defi}
The relative end depends on the choice of the prebalancing. We will denote the relative end as $\oint_{M\in \Mo} (S,\beta)$, or sometimes simply as $\oint_{M\in \Mo} S$, when the prebalancing $\beta$ is understood from the context.

\smallbreak
The \textit{ relative coend} of the pair $(S,\beta)$ is defined dually. This  is an object $C\in \Ac$ equipped with dinatural transformations $\pi_M:S(M,M)\xrightarrow{ . .} C$ such that
\begin{equation}\label{dinat:coen:module:left} 
\pentagonABxBCxCDxDExAE
{S(M,M)}
{S(\mathrm{id}_M,\mathrm{coev}_X\triangleright\mathrm{id}_M)}
{S(M,(X^*\otimes X)\triangleright M)}
{S(\mathrm{id}_M, m_{X,X^*,M})\vphantom{{{X^X}^X}^X}} 
{S(M,X\triangleright (X^*\triangleright M))}
{\beta^X_{M,X^*\triangleright M}}
{S(X^{*}\triangleright M,X^*\triangleright M)}
{\pi_{X^*\triangleright M}}
{C}
{\pi_M}
\end{equation}
for any $X\in \ca, M\in \Mo$, universal with this property. This means that, if $\widetilde{C}\in \Ac$ is another object with dinatural transformations  $\lambda_M: S(M,M) \xrightarrow{ . .} \widetilde{C}$ such that they satisfy \eqref{dinat:coen:module:left}, there exists a unique morphism $g:C\to \widetilde{C}$ such that $g \circ \pi_M=\lambda_M$. The relative coend will be denoted  $\oint^{M\in \Mo} (S,\beta)$, or simply as $\oint^{M\in \Mo} S$.

\medbreak

A similar definition can be made for \textit{right} $\ca$-module categories. Let $\Ac$ be a category, and $\No$ be a right $\ca$-module category endowed with a functor $S:\No^{\op}\times \No\to \Ac$ with a \textit{prebalancing}
$$ \gamma^X_{M,N}: S(M\triangleleft X, N)\to S(M, N\triangleleft {}^* X),$$
for any $M,N\in \No$, $X\in \ca$. 
\begin{defi} The \textit{relative end} for $S$ is an object $E\in \Ac$ equipped with dinatural transformations  $\lambda_N:E\xrightarrow{ .. } S(N,N)$ such that 
\begin{equation}\label{dinat:end:module:right}
\pentagonD
{S(N,N)}
{ S(\mathrm{id}_N,\mathrm{id}_N \triangleleft \,\mathrm{ev}_X)}
{S(N,N\triangleleft (X\otimes {^*}X))}
{S(\mathrm{id},\mathrm{id}\otimes m^{-1}_{N,X,{}^* X} )}
{S(N,(N\triangleleft X)\triangleleft {^*}X)}
{\gamma^X_{N,N\triangleleft X}}
{S(N\triangleleft X,N\triangleleft X)}
{\lambda_{N\triangleleft X}}
{E}
{\lambda_N}
\end{equation}
for any $N\in \No$, $X\in \ca$. We shall also denote this relative end by $\oint_{N\in \No} (S, \gamma)$.

Similarly, the \textit{relative coend} is an object $C\in \Bc$ with dinatural transformations $\lambda_N: S(N,N)\xrightarrow{ .. } C$  such that
\begin{equation}\label{dinat:coend:module:right}
\pentagonBAxBCxCDxDExAE
{S(N,N)}
{S(\mathrm{id}_N  \triangleleft\, \mathrm{coev}_X^{op},\mathrm{id}_N)}
{S(N  \triangleleft  ({^*}X\otimes  X),N)}
{S(m^{-1\,op}_{N, {^*}X, X}, \mathrm{id}_N)}
{S((N  \triangleleft  {}^*X)\triangleleft X,N )}
{\gamma^X_{N  \triangleleft  {^*}X, N}}
{S(N  \triangleleft  {^*}X,N  \triangleleft  {}^* X)}
{\lambda_{N  \triangleleft  {^*}X}}
{C}
{\lambda_N}
\end{equation}
for any $N\in \No$, $X\in \ca$. We shall also denote this relative coend by $\oint^{N\in \No} (S, \gamma)$.
\end{defi}

In the next Proposition we collect some results about the relative (co)end that will be useful.  The reader is referred to \cite[Prop. 3.3]{BM}, \cite[Prop. 4.2]{BM}.

\begin{prop}\label{compilat-r} Assume  $\Mo, \No$ are  left $\ca$-module categories, and $S, \widetilde{S}:\Mo^{\op}\times \Mo\to \Ac$ are functors equipped with $\ca$-prebalancings $$\beta^X_{M,N}: S(M,X\triangleright N)\to S(X^*\triangleright M,N),$$ $$ \widetilde{\beta}^X_{M,N}: \widetilde{S}(M,X\triangleright N)\to \widetilde{S}(X^*\triangleright M,N),$$ $X\in \ca, M,N\in \Mo$. The following assertions holds
\begin{itemize}

\item[(i)] Assume that the module ends $\oint_{M\in \Mo} (S,\beta), \oint_{M\in \Mo} (\widetilde{S},\widetilde{\beta})$ exist and have dinatural transformations $\pi, \widetilde{\pi}$, respectively. If  $\gamma:S\to \widetilde{S}$ is a natural transformation  such that 
\begin{equation}\label{gamma-ind}
 \widetilde{\beta}^X_{M,N} \gamma_{(M,X\triangleright N)}=\gamma_{(X^*\triangleright M,N)} \beta^X_{M,N},
\end{equation}
then there exists a unique map $\widehat{\gamma}: \oint_{M\in \Mo} (S,\beta)\to \oint_{M\in \Mo} (\widetilde{S},\widetilde{\beta})$ such that   $$\widetilde{\pi}_M \widehat{\gamma}= \gamma_{(M,M)} \pi_M$$ for any $M\in \Mo$. If $\gamma$ is a natural isomorphism, then  $\widehat{\gamma}$ is an isomorphism.

\item[(ii)] For any pair of $\ca$-module functors $(F, c), (G, d):\Mo\to \No,$ the functor $$\Hom_{\No}(F(-), G(-)): \Mo^{\op}\times \Mo\to \vect_\ku$$
has a canonical prebalancing  given by
\begin{equation}\label{beta-for-hom}
 \beta^X_{M,N}: \Hom_{\No}(F(M), G(X\triangleright N))\to \Hom_{\No}(F(X^*\triangleright M), G(N))
\end{equation}
$$ \beta^X_{M,N}(\alpha)= (ev_X\triangleright \id_{G(N)}) m^{-1}_{X^*,X,G(N)} (\id_{X^*}\triangleright d_{X,N}\alpha)c_{X^*,M},$$
for any $X\in \ca, M, N\in \Mo$.
There is an isomorphism $$ \Nat_{\!m}(F,G)\simeq \oint_{M\in \Mo} (\Hom_{\No}(F(-), G(-)), \beta).$$\qed
\end{itemize}
\end{prop}
The next result will be needed later. It follows from a combination of \cite[Prop. 3.3]{BM} (ii) and \cite[Lemma 3.6]{BM}.
Let $\Mo$ be a left $\ca$-module category, then $\No=\Mo^{\op}$ is a right $\ca$-module category with right action given by
$$M \triangleleft X= X^*\triangleright M.  $$
 Assume $\Ac$ is a category equipped with a functor
$$S:\No^{\op}\times \No\to \Ac,$$
together with a prebalancing
$$\beta^X_{M,N}: S(M \triangleleft X, N)\to S(M,N \triangleleft {}^*X).$$
The functor
$$\Hom_\Ac(S(-,-), U):\Mo^{\op}\times \Mo\to \vect_\ku, $$
has a natural prebalancing 
$$ \gamma^X_{M,N}:\Hom_\Ac(S(M, N\triangleleft {}^* X), U)\to \Hom_\Ac( S(M\triangleleft X, N),U),$$
$$\gamma^X_{M,N}(f)=f\circ \beta^X_{M,N}.$$
The proof of the next statement is similar to the proof of \cite[Prop. 3.3 (ii)]{BM}.
\begin{prop}\label{int-hom} If the coend  $\oint_{M\in \No} (S,\beta)$ exists, then for any object $U\in \Ac$,  the end $\oint_{M\in \Mo} \Hom_\Ac( S( -, -), U)$ exists, and there is an isomorphism
$$\oint_{M\in \Mo} \Hom_\Ac(S( -, -),U)\simeq  \Hom_\Ac( \oint^{M\in \No} (S,\beta), U).$$
Moreover, if $\oint_{M\in \Mo} \Hom_\Ac(S( -, -),U)$ exists for any $U\in \Ac$, then the coend $\oint^{M\in \No} (S,\beta)$ exists.\qed
\end{prop}

The proof of the next result is completely analogous to the proof in the case of usual coend.

\begin{lema}\label{right-exact-on-coend} Let $\Mo$ be a right $\ca$-module category, $S:\Mo^{\op}\times \Mo\to \ac$ a functor equipped with a prebalancing $\gamma$. If $F:\Ac\to \Ac'$ is a right exact functor, then there is an isomorphism
$$F(\oint^{M\in \Mo} (S,\gamma))\simeq \oint^{M\in \Mo} (F\circ S,F(\gamma)).$$\qed
\end{lema}

We shall also need the next result.
\begin{lema}\label{equivalence-on-coend} Let $\Mo, \No$ be right $\ca$-module categories and let $(J,c):\Mo\to\No$ be an equivalence of $\ca$-module categories. Assume that $S:\No^{\op}\times \No\to \ac$ is a functor equipped with a prebalancing $\beta$. The functor $S(J(-), J(-)):\Mo^{\op}\times \Mo\to \ac$ has a prebalancing $\gamma$ given by
$$\gamma^V_{M,N}=S(\id_{J(M)},c_{N,{}^*V})\beta^V_{J(M),J(N)}S(c^{-1}_{M,V},\id_{J(N)}) $$
for any $M,N\in \Mo$. There is an isomorphism
$$\oint^{M\in \Mo}  (S(J(-), J(-)),\gamma)\simeq \oint^{N\in \No}  (S,\beta).$$
\end{lema}
\pf We only sketch the proof. Let $\lambda_N:S(N,N)\to \oint^{N\in \No}  (S,\beta)$ the associated dinatural transformations. If we define $\pi_M=\lambda_{J(M)}$, for any $M\in \Mo$, one can verify that these are dinatural transformations and they satisfy \eqref{dinat:coend:module:right}. Universality of $\pi$ follows from the universality of $\lambda$ and the fact that $J$ is an equivalence of categories.
\epf

\subsection{Existence of (co)ends }\label{Section_coendExists}

In this Section we shall prove that certain (co)ends exists. All these (co)ends will be used in subsequent sections. We shall use ideas from \cite{Sh5}.

\medbreak

Let $\vc$ be a tensor category, $\Bc$  be  a right $\vc$-module category. Since our definition of module category includes that the action is right exact in each variable, hypothesis  of \cite[Thm. 2.24]{DSS} are fulfilled. This means that, there exists an algebra $A\in \vc$ such that 
\begin{equation}\label{eq-ma}
\Bc\simeq {}_{A}\vc
\end{equation} 
as module categories. Consider $\vc$ as a right $\vc$-module category with the regular action. We shall denote by $\Rex_\vc(\Bc, \vc)$ the category of right exact $\vc$-module functors. The functor 
$$\Phi:\Bc^{\op} \to \Rex_\vc(\Bc, \vc),$$
$$\Phi(B)(D)= \uhom(D,B)^*$$
is an equivalence of categories, since, under identification \eqref{eq-ma}, the functor $\Phi$ is the composition of equivalences
$$ ({}_A\vc)^{\op} \xrightarrow{( -)^*} \vc_A\xrightarrow{\, R \,} \Rex_\vc({}_{A}\vc,\vc).$$
Here $R:  \vc_{A}  \to \Rex_\vc({}_A\vc,\vc)$ is the functor given by $R(V)(W)=V\ot_{A} W.$
Let be $(F,c)\in \Rex_\vc(\Bc, \vc)$. In particular, we have natural isomorphisms
$$c_{B,X}: F(B \triangleleft X)\to F(B) \ot X,$$
for $B\in\Bc, X\in\vc$. Consider the functor
$$\Bc^{\op}\times \Bc\to \Bc, $$
$$(D, B)\mapsto D \triangleleft {}^*F(B). $$
 This functor posses a $\vc$-prebalancing
\begin{equation}\label{prebalancing-gamma}
    \gamma^V_{B,D}: D\triangleleft {}^*F(B \triangleleft V) \to D\triangleleft {}^*V\ot{}^*F(B ),
\end{equation} 
$$\gamma^V_{B,D}=\id_{D}\triangleleft   (\phi^l_{F(B),V})^{-1} {}^*(c^{-1}_{B,V}). $$
\begin{prop}\label{existence-coend0} Let $\Bc$  be  a right $\vc$-module category. For any right exact $\vc$-module functor $(F,c):\Bc\to \vc$ the coend
$$\widetilde{C}(F,c)=\oint^{B\in \Bc} (B \triangleleft {}^*F(B), \gamma)\in \Bc,$$
exists. Moreover, the functor $(F,c)\mapsto \widetilde{C}(F,c)$ is a quasi-inverse of $\Phi$.
\end{prop}
\pf For any $D\in \Bc$ we have
\begin{align*}
\Nat_{\! m}(\Phi(D), F)& \simeq \oint_{B\in \Bc} \Hom_\vc(\Phi(D)(B), F(B))\\
&\simeq \oint_{B\in \Bc} \Hom_\vc(\uhom_{\Bc}(B,D)^*, F(B))\\
&\simeq \oint_{B\in \Bc} \Hom_\vc( {}^*F(B), \uhom_{\Bc}(B,D))\\
&\simeq \oint_{B\in \Bc} \Hom_\Bc(B \triangleleft {}^*F(B), D)\\
&\simeq \Hom_\Bc( \oint^{B\in \Bc}  B \triangleleft {}^*F(B), D)\\
&=\Hom_{\Bc^{\op}}(D, \oint^{B\in \Bc}  B \triangleleft {}^*F(B)).
\end{align*} 
The first isomorphism is Proposition \ref{compilat-r} (ii), and the fifth isomorphism is Proposition \ref{int-hom}. Observe that, to prove the existence of the  fourth isomorphism, one has to check that the natural isomorphisms
$$\Hom_\vc( {}^*F(B), \uhom_{\Bc}(B,D)\simeq \Hom_\Bc(B \triangleleft {}^*F(B), D)$$
commute with the respective prebalancings and then use Proposition \ref{compilat-r} (i). This calculation is left to the reader.
This proves both, that the coend  $\widetilde{C}(F,c)$ exists and that the quasi-inverse of $\Phi$ is given by 
$$\bar{\Phi}:  \Rex_\vc(\Bc, \vc) \to\Bc^{\op}, $$ $$\bar{\Phi}(F,c)= \oint^{B\in \Bc} (B \triangleleft {}^*F(B), \gamma).$$\epf

For a pair of  module functors $(F,c), (\widetilde F, \tilde c):\Bc\to \vc$, the functor
$$\Bc^{\op}\times \Bc\to \vc, $$
$$ (A,B)\mapsto \widetilde{F}(B)\ot {}^*F(A)$$
 has a canonical prebalancing, given by
\begin{equation}\label{preb-s1}\begin{split} \beta^X_{A,B}: \widetilde{F}(B)\ot {}^*F(A \triangleleft X) \to \widetilde{F}(B \triangleleft  {}^*X)\ot  {}^*F(A),\\
\beta^X_{A,B}=(\tilde{c}^{-1}_{B, {}^*X}\ot \id)(\id_{\widetilde F(B)}\ot (\phi^l_{F(A),X})^{-1} {}^*(c^{-1}_{A,X})).
\end{split}
\end{equation}

As a consequence of Lemma \ref{right-exact-on-coend}, one can apply $\widetilde F $ to the coend in Proposition \ref{existence-coend0}. Thus we get the next result.
\begin{cor}\label{existence-coend1}  For any pair of right exact $\vc$-module functors $F, \widetilde F:\Bc\to \vc$, the coend
$$\oint^{B\in \Bc} (\widetilde F(B) \ot {}^*F(B), \beta)\in \Bc,$$
exists.\qed
\end{cor}

\section{ Hopf algebras in braided tensor categories}\label{Section_Hopf}

\subsection{Hopf algebras in braided tensor categories}
Let us briefly recall the notion of Hopf algebras in braided tensor categories, and how their corepresentation categories are again tensor categories. For more details, the reader is referred to \cite{Majid},  \cite{Majid2}, \cite{Ta} and references therein. 

\medbreak

A \textit{braided tensor category} is a pair $(\vc,\sigma)$ where $\vc$ is a tensor category and $\sigma_{V,W}:V\ot W\to W\ot V$ is a \textit{braiding}, that is, a family of natural isomorphisms satisfying
\begin{equation}\label{braiding1}\begin{split} \sigma_{V,U\ot W}=(\id_U\ot \sigma_{V,W})(\sigma_{V,U}\ot\id_W),\\
\sigma_{V\ot U, W}=(\sigma_{V,W}\ot\id_U) (\id_V\ot \sigma_{U,W}).
\end{split}
\end{equation}
\begin{rmk} Note that in the above axioms we are assuming that $\vc$ is a strict tensor category. 
\end{rmk}
The braiding fulfills the braid relation, whence the name. We illustrate this identity in terms of string diagrams, which we read bottom to top
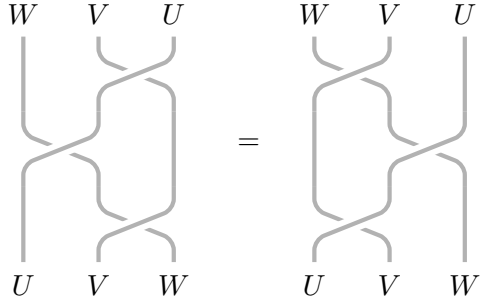
\begin{center}
        \begin{grform}
           \begin{scope}[scale = 1]
                \def\y{-0.3}\def\yn{-0.3+1}
                \draw (0 , \y) node {$U$};
                \draw (1 , \y) node {$V$};
                \draw (2 , \y) node {$W$};

                \def\y{0}\def\yn{1}
                \vLine{0}{\y}{0}{\yn}{\grau};
                \vLine{2}{\y}{1}{\yn}{\grau};
                \vLineO{1}{\y}{2}{\yn}{\grau};

                \def\y{1}\def\yn{2}
                \vLine{1}{\y}{0}{\yn}{\grau};
                \vLineO{0}{\y}{1}{\yn}{\grau};
                \vLine{2}{\y}{2}{\yn}{\grau};

                \def\y{2}\def\yn{3}
                \vLine{0}{\y}{0}{\yn}{\grau};
                \vLine{2}{\y}{1}{\yn}{\grau};
                \vLineO{1}{\y}{2}{\yn}{\grau};

                \def\y{3.3}\def\yn{4.3}
                \draw (0 , \y) node {$W$};
                \draw (1 , \y) node {$V$};
                \draw (2 , \y) node {$U$};
            \end{scope}  
        \end{grform}
                \quad=\quad
        \begin{grform}
            \begin{scope}[scale = 1]
                \def\y{-0.3}\def\yn{-0.3+1}
                \draw (0 , \y) node {$U$};
                \draw (1 , \y) node {$V$};
                \draw (2 , \y) node {$W$};

                \def\y{0}\def\yn{1}
                \vLine{1}{\y}{0}{\yn}{\grau};
                \vLineO{0}{\y}{1}{\yn}{\grau};
                \vLine{2}{\y}{2}{\yn}{\grau};

                \def\y{1}\def\yn{2}
                 \vLine{0}{\y}{0}{\yn}{\grau};
                \vLine{2}{\y}{1}{\yn}{\grau};
                \vLineO{1}{\y}{2}{\yn}{\grau};

                \def\y{2}\def\yn{3}
                \vLine{1}{\y}{0}{\yn}{\grau};
                \vLineO{0}{\y}{1}{\yn}{\grau};
                \vLine{2}{\y}{2}{\yn}{\grau};

                \def\y{3.3}\def\yn{4.3}
                \draw (0 , \y) node {$W$};
                \draw (1 , \y) node {$V$};
                \draw (2 , \y) node {$U$};
            \end{scope}
        \end{grform}
\end{center}

\begin{defi}
A \emph{bialgebra} in $\vc$ is a collection $(H,m,u,\Delta,\varepsilon)$, where $(H,m,u)$ is an algebra, $(H,\Delta,\varepsilon)$ is a coalgebra, and $\Delta, \varepsilon$ are algebra morphisms. I.e.
\begin{align}\label{definition:bialgebra} 
\Delta\circ m&=(m\ot m) (\id\ot \sigma_{H,H}\ot \id)  (\Delta\ot \Delta),\\
\Delta\circ u&= (u\ot u)\\
\varepsilon \circ m &= (\varepsilon\ot\varepsilon) \\
\varepsilon \circ u &= \id_\uno
\end{align}
We illustrate the first identity
\begin{center}
            \begin{grform}
                    \begin{scope}[scale = .8]
                       \dMult{1}{4}{2}{-1.5}{\grau}
                       \dMult{1}{0}{2}{1.5}{\grau}
                       \vLine{2}{1}{2}{3}{\grau}
                        
                        \draw (1 , -0.3) node {$H$};
                        \draw (1 , 4.3) node {$H$};
                        \draw (3 , -0.3) node {$H$};
                        \draw (3 , 4.3) node {$H$};
                    \end{scope}
                \end{grform}
                = \hspace{.2cm}
                \begin{grform}
                    \begin{scope}[scale = .8]
                       \dMult{0.5}{1.5}{1}{-1.5}{\grau}
                       \dMult{0.5}{2.5}{1}{1.5}{\grau}
                       \dMult{2.5}{1.5}{1}{-1.5}{\grau}{black}
                       \dMult{2.5}{2.5}{1}{1.5}{\grau}{black}
                       \vLine{0.5}{1.5}{0.5}{2.5}{\grau}
                       \vLine{2.5}{1.5}{1.5}{2.5}{\grau}
                       \vLineO{1.5}{1.5}{2.5}{2.5}{\grau}
                       \vLine{3.5}{1.5}{3.5}{2.5}{\grau}
                        
                        \draw (1 , -0.3) node {$H$};
                        \draw (1 , 4.3) node {$H$};
                        \draw (3 , -0.3) node {$H$};
                        \draw (3 , 4.3) node {$H$};
                   \end{scope}
                \end{grform}
\end{center}
    
\end{defi}
If $H$ is a bialgebra, then the space $\Hom_\vc(H,H)$ has a a convolution product with unit given by $u\circ \varepsilon$. If the identity $\id_H$ has an inverse $S$ under the convolution product, then $H$ is a \textit{Hopf algebra}, and $S$ is called the \textit{antipode}. The next Theorem is  well known, see for example \cite{Majid}.
\begin{teo}
If $H$ is a Hopf algebra in a braided tensor category $\vc$, then the category of left $H$-comodules ${}^H\vc$ is a  tensor category. Moreover, if $\vc$ is a finite tensor category, then ${}^H\vc$  is also  a finite tensor category, and the forgetful functor $\uc:{}^H\vc\to \vc$ is an exact faithful monoidal functor.\qed
\end{teo}

The tensor product of two left $H$-comodules is given as follows: If $(V,\rho_V)$, $(W, \rho_W)$ are objects in ${}^H\vc$, then the tensor product $V\ot W$ in $\vc$ has a left $H$-comodule structure given by
\begin{equation}\label{comodule-tensor-prod}\rho_{V\ot W}= (m\ot \id_{V\ot W})(\id_H\ot \sigma_{V,H}\ot \id_W)(\rho_V\ot \rho_W).
\end{equation}
\begin{center}
            \begin{grform}
                        \begin{scope}[scale = 1]
                                \dMult{-0.5}{2}{1}{1}{\grau}
                                \dAction{-0.5}{1}{2}{-1}{\grau}{black}
                                \dAction{0.5}{2}{2}{-1}{\grau}{black}
                                \vLine{-0.5}{1}{-0.5}{2}{\grau}
                                \vLineO{1.5}{0}{1.5}{3}{black}
                                \vLine{2.5}{0}{2.5}{3}{black}
                                \draw (0 , 3.3) node {$H$};
                                \draw (1.5 , -0.3) node {$V$};
                                \draw (2.5 , -0.3) node {$W$};
                                \draw (1.5 , 3.3) node {$V$};
                                \draw (2.5 , 3.3) node {$W$};
                        \end{scope}
                        
                \end{grform}
    \end{center}

If $(V, \rho_V)\in {}^H\vc$ then ${}^*V\in  {}^H\vc$. The coaction is given by 
\begin{equation}\label{antipode-com} 
\rho_{{}^*V}=\big((S\ot \id)\sigma_{{}^*V,H}\ot \ev_V\big)\big(\id_{{}^*V}\ot \rho_V\ot \id_{{}^*V}\big)\big(\coev_V\ot\id_{{}^*V}\big).
\end{equation}
The coaction of the right dual $V^*$ is defined similarly.

\medbreak

For any $V\in \vc$ we can endow $V$ with a trivial $H$-comodule structure, given by
$$\rho^t_V:V\to H\ot V, \quad \rho^t_V=u\ot \id_V.$$

\subsection{Yetter-Drinfeld modules}\label{Section_YD}

Let $H$ be a finite dimensional Hopf algebra in $\vect_\ku$. Then we shall denote by ${}^H_H\YD$ the category of finite-dimensional \textit{Yetter-Drinfeld modules.}  An object $V\in {}^H_H\YD(\vect_\ku)$ is a left $H$-module $\cdot:H\otk V\to V$, and a left $H$-comodule $\lambda:V\to H\otk V$ such that
\begin{equation}\label{yd-module} \lambda(h\cdot v)=h\_1 v\_{-1} \Sc(h\_3)\ot h\_2\cdot v\_0,
\end{equation}
for any $h\in H, v\in V$.  If $V\in {}^H_H\YD(\vect_\ku)$, the map
$\sigma_X: V\otk X\to X\otk V$, given by $\sigma_X(v\ot x)=v\_{-1}\cdot x\ot v\_0$ is a half-braiding for $V$.\\

This notion had been generalized in \cite{Besp95} to Hopf algebras $H$ inside a braided tensor category $\vc$:

\begin{defi}\label{def_YD}
Let $H$ be a Hopf algebra in a braided tensor category $\vc$. Then an $H$-$H$-\emph{Yetter-Drinfeld module} $V$ is an object $V\in\vc$, together with a structure of  $H$-module $m:H\otimes V\to V$  and  structure of  $H$-comodule in $\rho:V\to H\otimes V$ in the tensor category $\vc$, compatible in the following way
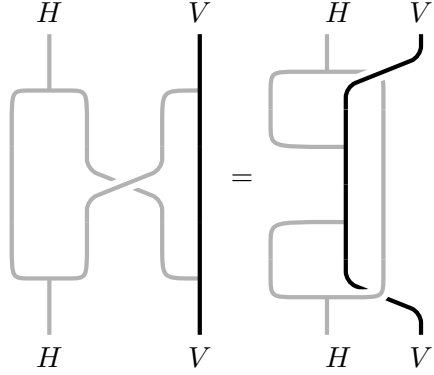
\begin{center}
                \begin{grform}
                    \begin{scope}[scale = 1]
                        \dMult{0.5}{1.5}{1}{-1.5}{\grau}
                        \dMult{0.5}{2.5}{1}{1.5}{\grau}
                        \dAction{2.5}{1.5}{0.5}{-1.5}{\grau}{black}
                        \dAction{2.5}{2.5}{0.5}{1.5}{\grau}{black}
                        \vLine{0.5}{1.5}{0.5}{2.5}{\grau}
                        \vLine{2.5}{1.5}{1.5}{2.5}{\grau}
                        \vLineO{1.5}{1.5}{2.5}{2.5}{\grau}
                        \vLine{3}{1.5}{3}{2.5}{black}
                        \draw (1 , -0.3) node {$H$};
                        \draw (1 , 4.3) node {$H$};
                        \draw (3 , -0.3) node {$V$};
                        \draw (3 , 4.3) node {$V$};
                    \end{scope}
                \end{grform}
                =
                \begin{grform}
                    \begin{scope}[scale = 1]
                       \vLine{2}{0}{1}{1}{black}
                       \dMultO{0}{1}{1.5}{-1}{\grau}
                       \dAction{0}{1}{1}{1}{\grau}{black}
                       \dAction{0}{3}{1}{-1}{\grau}{black}
                       \dMult{0}{3}{1.5}{1}{\grau}
                       \vLineO{2}{4}{1}{3}{black}
                       \vLine{1.5}{1}{1.5}{3}{\grau}
                        \draw (0.9 , -0.3) node {$H$};
                        \draw (0.9 , 4.3) node {$H$};
                        \draw (2 , -0.3) node {$V$};
                        \draw (2 , 4.3) node {$V$};
                    \end{scope} 
                \end{grform}
\end{center}
\end{defi}

The category ${^H_H}\YD(\vc)$ consists of  Yetter-Drinfeld modules and of $H$-linear and $H$-colinear morphisms. It becomes a tensor category with the usual tensor product of $H$-modules and $H$-comodules $V\otimes W$. An useful feature of this tensor category is that, it admits, by construction, a braiding
\begin{align*}
     c_{(V, m_V, \rho_V), (W, m_W, \rho_W)} 
    &: V \otimes W \rightarrow W \otimes V \\
\intertext{which is given on objects $(V, m_V, \rho_V)$ and $(W, m_W, \rho_W)$ by}
     c_{(V, m_V, \rho_V), (W, m_W, \rho_W)} 
    &:= (\rho_W \otimes \id_V ) \circ (\id_H \otimes
    c_{V,W} ) \circ (\rho_V \otimes \id_W )\,\, \\
\end{align*}
and is invertible if $H$ has a bijective antipode
\begin{align*}
(c_{(V, m_V, \rho_V), (W, m_W, \rho_W)} )^{-1} 
&:= c^{-1}_{V,W} 
\circ (m_W \otimes \id_V ) \circ (c^{-1}_{H,W} \otimes \id_V ) \\
&\circ (\id_W \otimes S^{-1} \otimes \id_V ) \circ (\id_W \otimes \rho_V ).
\end{align*}
 If $\vc$ is rigid, then the dual object in $\vc$ with the standard dual action and coaction  gives a dual object in ${^H_H}\YD(\vc)$. The structure is summarized in the following statement proven in \cite{Besp95}:
\begin{teo}
Let $H$ be a Hopf algebra in $\vc$. 
The Yetter-Drinfeld modules over $H$ in $\vc$ have
a natural structure of a
braided tensor category ${^H_H}\YD(\vc)$. If $\vc$ is rigid, then so is ${^H_H}\YD(\vc)$.
\end{teo}

Recall that the \emph{Drinfeld center}  
$\mathcal{Z}(\mathcal{B})$ is a braided tensor category associated to any tensor category $\mathcal{B}$.
The construction of Yetter-Drinfeld modules gives a realization of a special case of the construction of a \emph{relative Drinfeld center}  
$\mathcal{Z}_\vc(\mathcal{B})$, a braided tensor category associated to any tensor category $\mathcal{B}$ with a braided subcategory $\vc\hookrightarrow \mathcal{Z}(\mathcal{B})$. More precisely we have 
$${^H_H}\YD(\vc)=\mathcal{Z}_\vc({^H}\vc)$$

\section{Monoidal module categories}

Given $\vc$, $\Bc$ tensor categories. We shall  define what means that a tensor category $\vc$ \textit{acts} on $\Bc$.

\begin{defi}\label{definition:monoidal-mod} We call $\Bc$ a \emph{monoidal left $\vc$-module category} if $\Bc$ is a left  $\vc$-module category, with action given by $\triangleright:\vc\times \Bc\to \Bc, $ such that
 $\uno\triangleright\uno\simeq \uno,$
and there are natural isomorphisms
$l_{X,B}: X \triangleright B\xrightarrow{\;\simeq\; } (X\triangleright \uno)\ot B,$
for any $X\in \vc$, $B\in \Bc$, such that
\begin{equation}\label{tens-acting0}  l_{\uno, B}=\id_B, \quad l_{X,\uno}=\id_{X\triangleright \uno},
\end{equation}
\begin{equation}\label{tens-acting}  
            \begin{tikzpicture}[commutative diagrams/every diagram]
            \node (P1) at (180:5cm) 
            {${X\triangleright(Y \triangleright B)}$};
            \node (P2) at (180-51:5cm) 
            {${(X\triangleright\uno)\otimes (Y \triangleright B)}$};
            \node (P3) at (180-51*2:5cm) 
            {$(X\triangleright\uno)\otimes ((Y \triangleright \uno)\otimes B)$};
            \node (P4) at (180-51*3:5cm) 
            {$((X\triangleright \uno)\otimes (Y \triangleright \uno))\otimes B$};
            \node (P5) at (180-51*4:5cm) 
            {$(X\triangleright (Y \triangleright \uno))\otimes B$};
            \node (P6) at (180-51*5:5cm)
            {${((X\otimes Y) \triangleright \uno)\otimes B}$};
            \node (P7) at (180-51*6:5cm) 
            {$(X\otimes Y) \triangleright B$};
\path[commutative diagrams/.cd, every arrow, every label]
            (P1) edge node [swap] {$l_{X,Y\triangleright B}$} (P2)
            (P2) edge node [swap] {$\mathrm{id}_{X\triangleright \uno}\ot  l_{Y,B}$} (P3)
            (P3) edge node [swap] {\rotatebox{-31}{$\sim$}} (P4)
            (P1) edge node {$m_{X,Y,B}^{-1}$} (P7)
            (P7) edge node {$l_{X\ot Y,B}$} (P6)
            (P6) edge node {$m_{X,Y,\uno}\otimes \mathrm{id}_B$} (P5)
            (P5) edge node {$l_{X,Y\triangleright  \uno}\ot\id_B$} (P4); 
            \end{tikzpicture}
\end{equation}

for any $X, Y\in \vc$, $B\in \Bc$. Similarly, we shall say that $\Bc$ is a monoidal right $\vc$-module category, with action $\triangleleft:\Bc\times  \vc\to \Bc, $ such that
$$\uno\triangleleft \uno\simeq \uno,$$
and there are natural isomorphisms
$l_{B,X}: B \triangleleft X\to B\ot (\uno \triangleleft X),$
such that for any $X, Y\in \vc$, $B\in \Bc$
\begin{equation}\label{tens-acting0r}  l_{ B,\uno}=\id_B, \quad l_{\uno, X}=\id_{\uno \triangleleft X },
\end{equation}
\begin{equation}\label{tens-acting2}  
            \begin{tikzpicture}[commutative diagrams/every diagram]
            \node (P1) at (180:5cm) 
            {$(B\triangleleft X)\triangleleft Y$};
            \node (P2) at (180-51:5cm) 
            {$(B\triangleleft X)\otimes (\uno\triangleleft Y)$};
            \node (P3) at (180-51*2:5cm) 
            {$(B\otimes (\uno\triangleleft X))\otimes (\uno\triangleleft Y)$};
            \node (P4) at (180-51*3:5cm) 
            {$B\otimes ((\uno\triangleleft X)\otimes (\uno\triangleleft Y))$};
            \node (P5) at (180-51*4:5cm) 
            {$B\otimes ((\uno\triangleleft X)\triangleleft Y))$};
            \node (P6) at (180-51*5:5cm)
            {$B\otimes (\uno\triangleleft (X\otimes Y))$};
            \node (P7) at (180-51*6:5cm) 
            {$B\triangleleft (X\otimes Y)$};
\path[commutative diagrams/.cd, every arrow, every label]
            (P1) edge node [swap] {$l_{B\triangleleft X,Y}$} (P2)
            (P2) edge node [swap] {$l_{B,X}\ot \mathrm{id}_{\uno \triangleleft Y}$} (P3)
            (P3) edge node [swap] {\rotatebox{-31}{$\sim$}} (P4)
            (P1) edge node {$\tilde{m}_{B,X,Y}^{-1}$} (P7)
            (P7) edge node {$l_{B, X\ot Y}$} (P6)
            (P6) edge node {$\mathrm{id}_B\otimes\tilde{m}_{\uno,X,Y}$} (P5)
            (P5) edge node {$\id_B\ot l_{\uno \triangleleft X, Y}$} (P4); 
            \end{tikzpicture}
\end{equation}
\end{defi}

\begin{exa}\begin{itemize}
\item[(i)]  Any tensor category $\vc$ is a monoidal $\vect_\ku$-module category, with the canonical action
\begin{equation*}
 \bullet: \vect_\ku \times \vc\to \vc.
\end{equation*}

\item[(ii)] Let $\vc$ be a tensor category. Any tensor subcategory $\Do$ acts on  $\vc$ by
$$X \triangleright Y=X\ot Y.$$
In particular any tensor category acts on itself.
\item[(iii)] Let $(\vc, \sigma)$ be a braided tensor category and $\Do\subseteq \vc$ be a tensor subcategory. Then $\Do^{\rev}$ acts on  $\vc$ as
$$X \triangleright Y=Y\ot X.$$
In this case $l_{X,B}=\sigma_{B,X}$.

\item[(iv)] Let $\vc$ be a tensor category and $C\in\vc$ be a coalgebra. Then the category ${}^C\vc$, of left $C$-comodules in $\vc$, is a right $\vc$-module category. The action is given as follows. If $(W,\rho)\in {}^C\vc$, then $W \triangleleft Y=W\ot Y$, where the coaction on $W\ot Y$ is given by $\rho\ot \id_Y.$

\end{itemize} 
\end{exa}

\begin{defi} If $\Bc$, $\Bc'$ are right monoidal $\vc$-module categories,   a \textit{monoidal module functor} is a collection  $(F,c,\xi):\Bc\to \Bc'$ where $(F,c)$ is a $\ca$-module functor, $(F, \xi)$ is a monoidal functor and equation \begin{equation}\label{monoidal-mod-functor} l_{F(B),V}\, c_{B,V}= (\id_{F(B)}\ot c_{\uno, V}) \xi^{-1}_{B, \uno \triangleleft V} F(l_{B,V}),
\end{equation}
is fulfilled, for any $B\in \Bc$, $V\in \ca$.
\end{defi}

\subsection{Sections of monoidal functors}

\begin{defi}\label{def_section}  Let $\vc, \Bc$ be tensor categories. If $F:\Bc\to \vc$ is a tensor functor, a \textit{section} to $F$ is a right exact tensor functor $G:\vc\to \Bc$ such that $F\circ G\simeq \Id_\vc$ as monoidal functors. 
\end{defi} 
\begin{exa} \begin{enumerate}
\item Any (linear) tensor functor $F:\Bc \to \vect_\ku$ has a canonical section given by
$$G:\vect_\ku \to \Bc,$$
$$G(V)= V  \bullet \uno.$$
\medbreak
Here the action $\bullet:\vect_\ku\times \Bc\to \Bc$, is the one presented in  \eqref{vect-action}.

\item The previous example can be generalized to  other kinds of fiber functors. If $\vc$ is a tensor category that acts on another tensor category $\Bc$ then, any monoidal functor $F:\Bc\to \vc$ that is also a $\vc$-module functor has a section given by $G:\vc\to \Bc$, $G(V)=V\triangleright \uno$. See Proposition \ref{lm_section} below.

\medbreak

\item If $(\vc,\sigma)$ is a braided tensor category, and $H\in \vc$ is a Hopf algebra with unit given by $u:\uno\to H$, then the forgetful functor $\uc: {}^H\vc\to \vc$ has a section given by $G:\vc\to {}^H\vc$, $G(V)=(V,\rho^t_V)$. Here  $\rho^t_V=u\ot \id_V$ is the trivial comodule structure.

\medbreak
\item Let $H$ be a Hopf algebra and $R$ be a  Hopf algebra in the category of Yetter-Drinfeld modules ${}^H_H\YD$. Consider the corresponding Hopf algebra obtained by bosonization $R\# H$. Let us consider the functor $$F:\Rep(R\# H)\to \Rep(H), \quad F(V)=V.$$ The action of $H$ on $V$ is given by $h\cdot v= (1\# h)\cdot v$. The functor $F$ has a section given by $G:\Rep(H)\to \Rep(R\# H)$ given by $G(V)=V$, where the action of $R\# H$ on $V$ is given by
$$(r\# h)\cdot v=\varepsilon(r) h\cdot v.$$

\item Let $(\vc,\sigma)$ be a braided tensor category. The forgetful functor from the center of $\vc$, $\uc :Z(\vc)\to \vc$,  has a section given by the inclusion $\vc \hookrightarrow Z(\vc)$, $V\mapsto (V,\sigma)$.

\end{enumerate}
\end{exa}

\begin{prop}\label{lm_section} Let $\Bc, \vc$ be tensor categories, and $(F,\xi):\Bc\to \vc$ be a monoidal functor. The following notions are equivalent:
\begin{itemize}
\item[(i)] The functor $F$   has a section;

\item[(ii)] $\Bc$ is a right monoidal $\vc$-module category (in the sense of Definition \ref{definition:monoidal-mod}) and $F$ is a  monoidal module  functor.
\end{itemize}
\end{prop}
\pf Let us only give a sketch of the proof.  (i) implies (ii): Assume that $(G, \zeta):\vc\to \Bc$ is a monoidal section to $F$. Let $\alpha:F\circ G\to \Id$ be a natural monoidal isomorphism. Define the right action of $\vc$ on $\Bc$ as
$$B \triangleleft V:= B\ot G(V), $$
for any $B\in \Bc$, $V\in \vc$. The associativity of this action is given by
$$ m_{B,V,W}: B \triangleleft (V\ot W)\to (B \triangleleft  V) \triangleleft W,$$
$$m_{B,V,W}= \id_B\ot \zeta^{-1}_{V,W},$$
for any $B\in \Bc$, $V, W\in \ca$. Turns out that $\Bc$ is a monoidal module category with isomorphisms
$$l_{B,V}: B \triangleleft V\to B\ot (\uno \triangleleft V),$$
$$l_{B,V}=\id_{B\ot G(V)}.$$

With this action, $F$ is a module functor. The module structure of the functor $F$ is given by
$$ c_{B,V}:F(B \triangleleft V)\to F(B)\ot V,$$
$$ c_{B,V}= (\id_{F(B)}\ot \alpha_V) \xi^{-1}_{B,G(V)} F(l_{B,V}).$$
Since $c_{\uno,V}=\alpha_V$, then it follows that Equation \eqref{monoidal-mod-functor} is fulfilled, that is $(F,\xi, c)$ is a monoidal module functor. 

\medbreak
Let us prove now that (ii) implies (i): Assume that $(F,\xi, c):\Bc\to \ca$ is a monoidal $\ca$-module functor. Define 
$G:\vc\to \Bc$, $G(V)=\uno \triangleleft V $. Axiom \eqref{tens-acting2} implies that  isomorphisms
$$l^{-1}_{\uno\triangleleft V, W }: G(V)\ot G(W)  \to G(V\ot W),$$
endow $G$ with the structure of monoidal functor. Define natural isomorphisms $\alpha_V:F(G(V))=F(\uno \triangleleft V)\to V$, $\alpha_V=c_{\uno, V}$. Let us check that they are monoidal natural isomorphisms. For this, we need to verify that equation
\begin{equation}\label{nat-c1} (c_{\uno, V}\ot c_{\uno, W})= c_{\uno, V\ot W} F(l^{-1}_{\uno \triangleleft V,W}) \xi_{\uno \triangleleft V,\uno \triangleleft W}
\end{equation}
is satisfied. Using that $c$ satisfies diagram \eqref{modfunctor11}, it follows that the right hand side of Equation  \eqref{nat-c1}  is equal to
\begin{align*} &=(c_{\uno, V}\ot \id_W) c_{\uno \triangleleft V, W} F(l^{-1}_{\uno \triangleleft V,W}) \xi_{\uno \triangleleft V,\uno \triangleleft W}\\
&=(c_{\uno, V}\ot \id_W) (\id_{F(\uno \triangleleft V)}\ot c_{\uno, W} ) \xi^{-1}_{\uno \triangleleft V,\uno \triangleleft W} F(l_{\uno \triangleleft V,W}) F(l^{-1}_{\uno \triangleleft V,W}) \xi_{\uno \triangleleft V,\uno \triangleleft W}\\
&=(c_{\uno, V}\ot c_{\uno, W}).
\end{align*}
The second equation follows from \eqref{monoidal-mod-functor}
\epf

\begin{rmk} Constructions in Proposition \ref{lm_section} are reciprocal in the following sense. If $\Bc$ is a right monoidal $\ca$-module category and $G:\ca\to \Bc$ is the monoidal functor $G(V)= \uno\triangleleft V$, then one can endow $\Bc$ with a  right $\ca$-module structure
$$ B\blacktriangleleft V= B\ot G(V),$$
for any $B\in \Bc$, $V\in \ca$. Turns out, that the identity functor $\Id: (\Bc, \triangleleft)\to (\Bc, \blacktriangleleft)$ is an equivalence of module categories.
\end{rmk}

\subsection{Some natural module transformations}\label{subsection:nat-bimod} It is known that natural transformations between additive functors are additive. For example, if $(F,\xi)$ is an additive monoidal functor, then, natural isomorphisms
$$\xi_{A,B}:F(A)\ot F(B)\to F(A\ot B), $$
are additive in each variable. In this section, we shall study what happens when $F$ is a monoidal module functor. We aim at proving that $\xi$ is a \textit{module} natural transformation in  the second variable. 

\begin{lema}\label{left-mod-nat} If $\Bc$ is a right monoidal $\ca$-module category, then for any $A\in \Bc$, the functor $L_A:\Bc\to \Bc$, $L_A(B)=A\ot B$ is a $\ca$-module functor.
\end{lema}
\pf For any $B\in \Bc, V\in \ca$, the module structure  of the functor $L_A$ is given by
$$\eta_{B,V}: A\ot (B \triangleleft V)\to (A\ot B) \triangleleft V,$$
$$\eta_{B,V}= l^{-1}_{A\ot B, V}(\id_A\ot l_{B,V})$$ 
Let us prove that $\eta$ satisfies   \eqref{modfunctor11}. Using that isomorphisms $l$ satisfy diagram \eqref{tens-acting2}, it follows that for any $C\in \Bc$ and any $V, W\in \ca$ we have that
\begin{equation}\label{on-l-11}
l^{-1}_{C, V\ot W}= l^{-1}_{C\triangleleft V,  W}(l^{-1}_{C,V}\ot \id_{\uno \triangleleft W})(\id_C\ot l_{\uno \triangleleft V,W}).
\end{equation}
 On one hand we have that 
\begin{align*}
\eta_{B,V\ot W}=& l^{-1}_{A\ot B, V\ot W}(\id_A\ot l_{B,V\ot W})\\
=& l^{-1}_{(A\ot B)\triangleleft V,  W}(l^{-1}_{A\ot B,V}\ot \id_{\uno \triangleleft W})(\id_{A\ot B}\ot l_{\uno \triangleleft V,W}) (\id_{A\ot B}\ot l^{-1}_{\uno \triangleleft V,W})\\
& (\id_A\ot l_{B,V}\ot \id_{\uno \triangleleft W})(\id_A\ot l_{B\triangleleft V,W})\\
=&l^{-1}_{(A\ot B)\triangleleft V,  W}(l^{-1}_{A\ot B,V}\ot \id_{\uno \triangleleft W})  (\id_A\ot l_{B,V}\ot \id_{\uno \triangleleft W})(\id_A\ot l_{B\triangleleft V,W})
\end{align*}
The second equation follows from using \eqref{on-l-11}. On the other hand $$(\eta_{B,V}\triangleleft \id_W) \eta_{B\triangleleft V, W}$$ is equal to
\begin{align*} =&(l^{-1}_{A\ot B, V} \triangleleft\id_W)(\id_A\ot l_{B,V}\ot \id_{\uno \triangleleft W}) l^{-1}_{A\ot (B \triangleleft V), W} (\id_A\ot l_{B \triangleleft V, W})\\
=&(l^{-1}_{A\ot B, V} \triangleleft\id_W)  l^{-1}_{A\ot B \ot (\uno\triangleleft V), W}  (\id_A\ot l_{B,V}\ot \id_{\uno \triangleleft W})  (\id_A\ot l_{B \triangleleft V, W})\\
=&l^{-1}_{(A\ot B)\triangleleft V,  W}(l^{-1}_{A\ot B,V}\ot \id_{\uno \triangleleft W})  (\id_A\ot l_{B,V}\ot \id_{\uno \triangleleft W})(\id_A\ot l_{B\triangleleft V,W}).
\end{align*}
The second and third equalities follow from the naturality of $l$.
\epf

\begin{prop}\label{bimodule-funct} Assume that $\Bc$ is a monoidal right $\ca$-module category and $(F,  c, \xi):\Bc\to \ca$ a monoidal module functor. The following assertions hold.
\begin{enumerate}
\item[(i)] The category $\Bc$ has structure of $\ca$-bimodule category.

\item[(ii)]  There are $\ca$-module funtors
$$H, \widetilde{H}:\Bc\boxtimes_\ca \Bc\to \ca,$$ such that
$$ H(A\boxtimes B)=F(A\ot B), \quad \widetilde{H}(A\boxtimes B)=F(A)\ot F(B),$$
for any $A, B\in \Bc$. 
\item[(iii)] The monoidal structure of the functor $F$ defines a natural module isomorphism $\xi:\widetilde{H} \to H$.
\end{enumerate}
\end{prop}
\pf (i). From Proposition \ref{lm_section} there is a section $G:\ca\to \Bc$ of the functor $F$. Define the left action $\triangleright:\ca\times \Bc\to \Bc$ as $V\triangleright B=G(V)\ot B$. With this action $\Bc$ becomes a $\ca$-bimodule category. Thus, we can consider the Deligne tensor product over $\ca$, $\Bc\boxtimes_\ca \Bc$. This category has an obvious right $\ca$-module category structure.

One can prove that functors $\Bc\times \Bc\to \ca$ given by
$$(A,B)\mapsto F(A\ot B), \quad (A,B)\mapsto F(A)\ot F(B),$$
are $\ca$-balanced, thus defining functors $H$ and $\widetilde{H}$. The (right) $\ca$-module structure of the functor $\widetilde{H}$ is the one inherited from the functor $F$. The module structure of the functor $H$ is more involved. One can prove that natural transformations
$$e_{A,B,V}: F(A\ot (B\triangleleft V))\to F(A\ot B) \ot V,$$
$$e_{A,B,V}= c_{A\ot B, V} F(l^{-1}_{A\ot B, V}(\id_A\ot l_{B,V})) $$
are $\ca$-balanced, thus defining natural transformations 
$$e_{X,V}: H(X\triangleleft V)\to H(X)\triangleleft V,$$
for any $X\in \Bc\boxtimes_\ca \Bc$, $V\in \ca$. One can prove also that $(H,e)$ is a module functor. Note that $e$ is the module structure resulting from the composition of module structures of $F$ and the functor $L_A$, presented in Lemma \ref{left-mod-nat}. This proves (ii).

(iii). To prove that $\xi$ is a natural module transformation in the second variable, we need to verify that the diagram
\begin{equation}\label{coproduct-cf}
\xymatrix{
 \widetilde{H}(A,B\triangleleft V) \ar[d]_{\id_{F(A)}\ot c_{B,V}}\ar[rr]^{\xi_{A,B\triangleleft V}}&& H(A ,B\triangleleft V) 
 \ar[d]^{e_{A,B,V}} \\
 \widetilde{H}(A,B) \ot V \ar[rr]_{\xi_{A,B}\ot \id_V }&& H(A ,B) \ot V,}
\end{equation}
is commutative for any $A, B\in \Bc$, $V\in \ca$. Here $$e_{A,B,V}= c_{A\ot B, V} F(l^{-1}_{A\ot B, V}(\id_A\ot l_{B,V})). $$ 
We have
\begin{align*} e_{A,B,V} &\xi_{A,B\triangleleft V}=  c_{A\ot B, V} F(l^{-1}_{A\ot B, V}(\id_A\ot l_{B,V}))\xi_{A,B\triangleleft V}\\
&=(\id\ot c_{\uno, V}) \xi^{-1}_{A\ot B, \uno\triangleleft V} F(\id_A\ot l_{B,V}) \xi_{A,B\triangleleft V}\\
&=(\id\ot c_{\uno, V}) \xi^{-1}_{A\ot B, \uno\triangleleft V} \xi_{A,B\ot (\uno \triangleleft V)}(\id_{F(A)}\ot F(l_{B,V}))\\
&=(\id\ot c_{\uno, V}) (\xi_{A,B}\ot \id_{F(\uno \triangleleft V)})(\id_{F(A)}\ot \xi^{-1}_{B, \uno \triangleleft V}) (\id_{F(A)}\ot F(l_{B,V}))\\
&=(\xi_{A,B}\ot \id_{ V}) (\id_{F(A)}\ot (\id_{F(B)}\ot c_{\uno, V})\xi^{-1}_{B, \uno \triangleleft V}  F(l_{B,V})) \\
&=(\xi_{A,B}\ot \id_{ V})  (\id_{F(A)}\ot c_{B, V}).
\end{align*}
The second equality follows from \eqref{monoidal-mod-functor}, the third equality follows from the naturality of $\xi$, the fourth one follows from \eqref{monoidal-functor}. The last equality follows  from \eqref{monoidal-mod-functor}.
\epf

\section{Fiber functors and reconstruction results}

It is a classical result that, out of a fiber functor, which is a monoidal, exact and faithful functor $F:\Bc\to \vect_\ku$, one can reconstruct a (usual) Hopf algebra $H$ such that $\Bc$ is monoidally equivalent to the category of left $H$-comodules.   We shall generalize these reconstruction theorems for a fiber functor $F:\Bc\to \vc$, where $\vc$ is an arbitrary braided tensor category acting monoidally on $\Bc$. The reconstruction of the Hopf algebra from these data will be described in a similar way as the case $\vc=\vect_\ku$.  See for example \cite{Majid}, \cite{PS}. The main new ingredient will be the use of the relative coend, introduced in \cite{BM}, and the splitting condition; that is, the existence of a section $G:\vc\to \Bc$ of $F$. While dealing with the relative (co)end, the primary new challenge will be demonstrate that some dinatural transformations satisfy Equation \eqref{dinat:coend:module:right}.

\subsection{ Coalgebras constructed from a coend}

Let $\vc$ be a  tensor category, and $\Bc$ be a right $\vc$-module category. For the rest of this section, $(F,c):\Bc\to \vc$ will denote  a right exact module functor. 

\medbreak

 Consider the functor
$$S^F:\Bc^{\op}\times \Bc\to \vc, $$
$$S^F(A,B)=F(B)\ot {}^*F(A).$$
This functor has a canonical prebalancing, given by
\begin{equation}\label{preb-s2}\begin{split} \beta^X_{A,B}: F(B)\ot {}^*F(A \triangleleft X) \to F(B \triangleleft  {}^*X)\ot  {}^*F(A),\\
\beta^X_{A,B}=(c^{-1}_{B, {}^*X}\ot \id_{{}^*F(A)})(\id_{F(B)}\ot (\phi^l_{F(A),X})^{-1} {}^*(c^{-1}_{A,X})).
\end{split}
\end{equation}
We define
$$C(F,c)=C(F)=\oint^{B\in \Bc} (S^F,\beta)= \oint^{B\in \Bc} F(B)\ot {}^*F(B).$$
Let $\pi_B:F(B)\ot {}^*F(B)\to C(F)$ be the associated dinatural transformations. Existence of this coend follows from Corollary \ref{existence-coend1}.
The following Lemma is a generalization of \cite[Lemma 2.1.9]{PS}. Recall that, for any $M, N\in \Bc$, $X\in \vc$ we have natural isomorphisms
$$\psi^V_{F(M),F(N)}:\Hom_\vc(F(M)\ot {}^*F(N),V)\to \Hom_\vc(F(M),V\ot F(N))$$
described in \eqref{duals1}.
\begin{lema}\label{a-omega} The functor 
$V\mapsto \Nat_{\! m}(F, V\ot F)$
is represented by $C(F)$. That is, for any $V\in\vc$ there are natural isomorphisms
\begin{equation}\label{omega-isom}
\begin{split} \omega: \Hom_\vc(C(F),V)\to \Nat_{\! m}(F, V\ot F),\\
\omega(g)_B= (g\ot \id_{F(B)})\psi^{C(F)}_{F(B),F(B)}(\pi_B).
\end{split}
\end{equation}  
\end{lema}
\pf Fix some object $V\in \vc$. The functor $V\ot F$ is a right module functor, then we can consider the prebalancing of the functor
$$\Hom_\vc(F( -),V\ot F(- )): \Bc\times \Bc^{\op}\to \vect_\ku, $$
given by
$$\gamma^X_{B,N} :\Hom_\vc(F(N\triangleleft {}^*X),V\ot F(B  ))\to   \Hom_\vc(F(N ),V\ot F(B\triangleleft X)),$$
$$\gamma^X_{N,B}(f)=(\id_V\ot c_{B, X})(f\ot \id_X) (c^{-1}_{N, {}^*X}\ot \id_X)(\id_{F(N)}\ot \coev_X).$$
Note that here we are considering $\B^{\op}$ as a left $\vc$-module category, with action $X\triangleright B:=B\triangleleft {}^*X$.
Proposition \ref{compilat-r} (ii) tells us that there is an isomorphism
$$\oint_{B\in \Bc}  (\Hom_\vc(F(B),V\ot F(B)), \gamma) \simeq \Nat_{\! m}(F, V\ot F). $$
From the proof of Proposition \ref{int-hom}, one can see that the dinatural transformations of the end $\oint_{B\in \Bc}  \Hom_\vc(F(B)\ot {}^*F(B),V)$ are given by
$$\widehat{\pi}_B:\oint_{B\in \Bc}  \Hom_\vc(F(B)\ot {}^*F(B),V)\to \Hom_\vc(F(B)\ot {}^*F(B),V), $$
$$\widehat{\pi}_B(g)= g\circ \pi_B.$$
Here, we are using identification $$\oint_{B\in \Bc}  \Hom_\vc(F(B)\ot {}^*F(B),V)= \Hom_\vc(\oint^{B\in \Bc}  F(B)\ot {}^*F(B),V).$$

 It follows by a straightforward computation that
$$\gamma^X_{M,N} (\psi^V_{F(M\triangleleft {}^*X ),F(N )}(f))=\psi^V_{F(M ),F(N\triangleleft X)}(f\circ \beta^X_{N,M} ),$$
for any $M, N\in \Bc$ and any $f\in \Hom_\Bc(F(M \triangleleft {}^*X)\ot {}^*F(N ),V)$. This implies, using Proposition \ref{compilat-r} (i) that, for any $V\in \vc$, there exists an isomorphism 
$$ \widehat{\psi}^V:\oint_{B\in \Bc}  \Hom_\vc(F(B)\ot {}^*F(B),V)\to \oint_{B\in \Bc}  \Hom_\vc(F(B),V\ot F(B)),$$
such that $\mu_B\circ\widehat{\psi}^V= \psi^V_{F(B),F(B)}\circ \widehat{\pi}_B.$ For any $V\in \vc$ we have that 
\begin{align*}  \Hom_\vc(C(F),V)&\simeq \oint_{B\in \Bc}  \Hom_\vc(F(B)\ot {}^*F(B),V)\\
&\simeq \oint_{B\in \Bc}  \Hom_\vc(F(B),V\ot F(B))\\
&\simeq \Nat_{\! m}(F, V\ot F).
\end{align*}
The first isomorphism is Proposition \ref{int-hom}. Tracing this chain of isomorphisms, one can see that, the composition coincides with $\omega$ defined by Equation \eqref{omega-isom}.
\epf

\begin{prop}\label{coalg-coend} Let $\vc$ be a  tensor category, and $\Bc$ be a right $\vc$-module category. Let $(F,c):\Bc\to \vc$ be  a right exact module functor. The object $C(F)\in\vc$ has a coalgebra structure $\Delta:C(F)\to C(F)\ot C(F)$, $\varepsilon:C(F)\to \uno$ determined by diagramms
\begin{equation}\label{coproduct-c}
\xymatrix{
 F(B)\ot {}^*F(B) \ar[d]_{\id\ot \mathrm{coev}_{F(B)} \ot\mathrm{id}}\ar[rr]^{\pi_B}&& C(F)
 \ar[d]^{\Delta} \\
 F(B)\ot {}^*F(B)\ot  F(B)\ot {}^*F(B) \ar[rr]_{\qquad\pi_B\ot \pi_B }&& C(F)\ot C(F),}
\end{equation}
\begin{equation}\label{counit-c}
\xymatrix{ F(B)\ot {}^*F(B)    \ar[rr]^{\pi_B}
\ar[dr]_{ \mathrm{ev}_{F(B)}  }&& C(F)
\ar[dl]^{\varepsilon }\\ & \uno,&}
\end{equation}
for any $B\in \Bc$.
\end{prop}
\pf It follows by a straightforward computation that, the maps
$$ \ev_{F(B)}:F(B)\ot {}^*F(B)\to \uno,  $$
$$(\pi_B\ot \pi_B )(\id\ot \coev_{F(B)} \ot\id): F(B)\ot {}^*F(B)\to C(F)\ot C(F), $$
are dinatural maps. It follows from Diagram \eqref{evaluation-mod} that $ \ev_{F(B)}$  satisfies \eqref{dinat:coend:module:right}.  It follows from Diagram \eqref{coevaluation-mod} and the fact that $\pi$ satisfies \eqref{dinat:coend:module:right} that, $(\pi_B\ot \pi_B )(\id\ot \coev_{F(B)} \ot\id)$ also satisfies Equation \eqref{dinat:coend:module:right}. Whence, the existence of $\Delta$ and $\varepsilon$ follow from the universal property of the relative coend. 
\medbreak

The proof that $(C(F),\Delta, \varepsilon)$ is a coalgebra is standard, and it follows from Diagrams \eqref{coproduct-c}, \eqref{counit-c}.
\epf 

Let $\vc$ be a tensor category and $C\in \vc$ be a coalgebra. The category ${}^C\vc$, of left $C$-comodules in $\vc$, is a right $\vc$-module category. The action is given as follows. If $(W,\rho)\in {}^C\vc$, then $W \triangleleft Y=W\ot Y$, where the coaction on $W\ot Y$ is given by $\rho\ot \id_Y.$
The forgetful functor $\uc:{}^C\vc\to \vc$ is a $\vc$-module functor. The next result says that the coalgebra reconstructed from the forgetful functor $\uc$ coincides with $C$.

\bigbreak

In this case, the functor $S^\uc$ has a prebalancing given by
$$\beta^X_{A,B}=\id_B\ot (\phi^l_{A,X})^{-1},$$
for any $A, B\in {}^C\vc$, $X\in \vc$. 
\begin{lema}\label{cc} There exists an isomorphism of coalgebras $C(\uc)\simeq C$ in $\vc$.
\end{lema}
\pf  In this particular case, if $T\in\vc$ is an object, and $\lambda_W:W\ot {}^*W\to T$ is a dinatural transformation, Equation \eqref{dinat:coend:module:right} writes as
\begin{equation}\label{new.3.5} 
\lambda_{W\ot {}^*V}(\id_W\ot (\phi^l_{W\ot {}^*V,V})^{-1}) =\lambda_W (\id_W\ot {}^*(\id_W\ot \coev_V)),
\end{equation}
for any $W\in {}^C\vc, V\in \vc.$

If $W\in {}^C\vc$ has comodule structure given by $\rho_W\to C\ot W$, then we define
$\pi_W:W\ot {}^*W\to  C,$ the morphisms given by $\pi_W=\bar{\psi}^C_{W,W}(\rho)$. Recall that, natural isomorphisms $\bar{\psi}$ were defined in \eqref{duals2}. Maps $\pi_W$ are dinatural transformations, and they satisfy Equation \eqref{new.3.5}.

Let $A\in \vc$ be another object, equipped with dinatural transformations $\lambda_W:W\ot {}^*W\to  A$ such that they satisfy Equation \eqref{new.3.5}. Define $h$ as the composition
$$h: C\xrightarrow{\id_C\ot {}^*\varepsilon}C\ot {}^*C   \xrightarrow{\lambda_C} A.$$
Since $\lambda$ is dinatural, and it satisfies  \eqref{new.3.5}, one can see that $h\circ \pi_W= \lambda_W$. This proves that, indeed, $C(\uc)\simeq C$.
\epf

\begin{prop}\label{factorization-F} Let $\Bc$ be a right $\vc$-module category, and $(F,c):\Bc\to \vc$ be a right exact module functor. Set $C=C(F)$ the coalgebra defined in Proposition \ref{coalg-coend}.
\begin{itemize}

\item[1.] The functor $(F,c):\Bc\to \vc$ factorizes as
\begin{equation}
\xymatrix{
&{}^C\vc\ar[dr]^-{\uc}\ar@{}[d]|{}&\\
\Bc\ar[ru]^-{\widehat{F}}\ar[rr]_-{F}&&\vc.
}
\end{equation} 
Here $\uc:{}^C\vc \to\vc$ is the forgetful functor, and the functor $\widehat{F}$ is a $\vc$-module functor.
\item[2.] If $F$ is exact (respect. faithful) then $\widehat{F}$ is exact (respect. faithful). 
\end{itemize}
\end{prop}
\pf 1. For any $A\in \Bc$, define $\widehat{F}(A)=(F(A),\rho_A),$ where 
$$\rho_A=\psi^C_{F(A),F(A)}(\pi_A):F(A)\to C\ot F(A).$$ Using diagrams \eqref{coproduct-c}, \eqref{counit-c}, one can prove that
$$(\id_C\ot \rho_A)\rho_A=(\Delta\ot \id)\rho_A,$$
$$(\varepsilon\ot\id)\rho_A=\id_{F(A)}.$$

\bigbreak

2. Clearly, if $F$ is faithful, then $\widehat{F}$ is faithful. Assume that $F$ is left exact. Let $f:A\to B$ be a morphism in $\Bc$ with kernel $k=\ker f: K\to A$. Since $F$ is left exact, then, $F(k)=\ker F(f)$. We only need to prove that, the morphism $F(k):F(K)\to F(A)$ is a map of $C$-comodules. That is, we need to prove the equality
\begin{equation}\label{exact-F} (\id\ot F(k))  \psi^C_{F(K),F(K)}(\pi_K)= \psi^C_{F(A),F(A)}(\pi_A)  F(k).
\end{equation}
Using the naturality of $\psi$, it follows that
\begin{align*}  \psi^C_{F(A),F(A)}(\pi_A)  F(k)&=  \psi^C_{F(K),F(A)}(\pi_A (F(k)\ot \id))\\
&=\psi^C_{F(K),F(A)}(\pi_K (\id \ot {}^*F(k)))\\
&=(\id \ot F(k)) \psi^C_{F(K),F(K)}(\pi_K ).
\end{align*}
The second equality follows from the dinaturality of $\pi$, the third equality follows again by the naturality of $\psi$. The proof that, if $F$ is right exact, then $\widehat{F}$ is right exact, follows similarly.
\epf

\subsection{Some auxiliary results} In this section we shall present some technical results that will be used later. In particular, these results will be used in the proof of our main result Theorem \ref{recons-th}.

\medbreak

Let  $\Bc$ be a right $\vc$-module category, and $(F,c):\Bc\to \vc$ be a module functor. To make calculations easier, we shall assume that the associativity of the $\vc$-module category $\Bc$ is trivial. Recall the relative coend
\begin{equation}\label{tilde-c} \widetilde{C}(F,c)=\widetilde{C}=\oint^{B\in \Bc} B \triangleleft {}^*F(B)\in \Bc,
\end{equation}
presented in Proposition \ref{existence-coend0}. Let us denote by
$$\lambda_B:B \triangleleft {}^*F(B) \to \widetilde{C}(F,c),$$
the associated dinatural transformations. Analogous to the definition of $\Delta$, given in 
\eqref{coproduct-c}, we shall define a map
\begin{equation}\label{coproduct-tilde-a}\delta: \widetilde{C}(F,c)\to \widetilde{C}(F,c)\triangleleft C(F,c).\end{equation}
\begin{lema} There exists a unique morphism $\delta: \widetilde{C}(F,c)\to \widetilde{C}(F,c)\triangleleft C(F,c)$ such that
\begin{equation}\label{coproduct-tilde-c}
\xymatrix{
 B\triangleleft {}^*F(B) \ar[d]^{\id\triangleleft(\coev_{F(B)} \ot\id)}\ar[rr]^{\lambda_B}&& \widetilde{C}
 \ar[d]^{\delta} \\
 B\triangleleft {}^*F(B)\ot  F(B)\ot {}^*F(B) \ar[rr]^{\qquad \,\,\,\,\,\,\lambda_B\triangleleft\,\pi_B }&& \widetilde{C}\triangleleft C.}
\end{equation}
\end{lema} 
\pf Let us denote $ d_{B,V}$ the following natural isomorphism
$$d_{B,V}:{}^*F(B\triangleleft V)\to  {}^*V\ot{}^*F(B) , $$
$$ d_{B,V}= \big(\phi^l_{F(B), V} \big)^{-1} {}^*(c^{-1}_{B,V}).$$
Recall from Equation \eqref{prebalancing-gamma}, that the prebalancing used in the coend $\widetilde{C}(F,c)$ is given by
$$\gamma^V_{B,D}: D\triangleleft {}^*F(B \triangleleft V) \to D\triangleleft {}^*V\ot{}^*F(B ), $$
$$\gamma^V_{B,D}=\id_{D}\triangleleft   d_{B,V}. $$
Using that $c$ satisfies \eqref{modfunctor11}, it is not difficult to prove that
\begin{equation}\label{about-d0} 
d_{B,{}^*V\ot V }=(\phi^l_{{}^*V,V}\ot \id)(\id_{{}^*V}\ot d_{B,{}^*V}) d_{B\triangleleft {}^*V, V },
\end{equation}
for any $B\in \Bc$, $V\in \vc$.
The naturality of $\phi^l$ and $c$ implies that
\begin{equation}\label{nat-phi90} 
F({}^*(\id_B\triangleleft \coev_V ))=({}^*(\coev_V)\ot \id_{{}^*F(B)}) d_{B,{}^*V\ot V },
\end{equation}
for any $B\in \Bc$, $V\in \vc$. Whence, using \eqref{coev-of-dual}, it follows that
\begin{equation}\label{f11} 
F({}^*(\id_B\triangleleft \coev_V )) d^{-1}_{B\triangleleft {}^*V, V }= (\ev_{{}^*V})\ot \id_{{}^*F(B)}) (\id_{{}^*V}\ot d_{B,{}^*V}).
\end{equation}
In this case, if $A\in \Bc$, a dinatural transformation $\nu_B: B\triangleleft {}^*F(B)\to A$ satisfies Equation \eqref{dinat:coend:module:right} if and only if
\begin{equation}\label{3.5-incase}  \nu_B \big( \id_B \triangleleft F({}^*(\id_B\triangleleft \coev_V )\big)= \nu_{B \triangleleft {}^*V}  \gamma^V_{B\triangleleft {}^*V,B}.
\end{equation}
Using that dinatural transformations $\lambda$ and $\pi$ both satisfy \eqref{dinat:coend:module:right} and using \eqref{f11}, it follows that dinatural transformation
$$(\lambda_B\triangleleft\,\pi_B) (\id_B\triangleleft(\coev_{F(B)} \ot\id_{{}^*F(B)}) $$
also satisfies \eqref{3.5-incase}. Now, existence of the map $\delta: \widetilde{C}(F,c)\to \widetilde{C}(F,c)\triangleleft C(F,c)$ satisfying \eqref{coproduct-tilde-c} follows from the universal property of the relative coend.
\epf

For any $A, B\in \Bc$ and $V\in \vc$, define
\begin{align}\label{duals01}\begin{split} \widetilde{\psi}^V_{B,A}:\Hom_\Bc(B\triangleleft {}^*V,A)\to \Hom_\Bc(B,A \triangleleft  V),\\
\widetilde{\psi}^V_{B,A}(f)=(f\triangleleft \id_V)(\id_B  \triangleleft \coev_V).
\end{split}
\end{align}
Compare with natural isomorphisms given in \eqref{duals1}. For any $B\in \Bc$, define also
\begin{equation}\label{coact1} \widetilde{\rho}_B= \widetilde{\psi}^{F(B)}_{B,\widetilde{C}}(\lambda_B):B\to \widetilde{C} \triangleleft F(B).
\end{equation}

One could check that, using Lemma \ref{right-exact-on-coend}, whenever $F$ is right exact, $F(\widetilde{C})=C$. The dinatural transformations of $F(\widetilde{C})$ is given by $F(\lambda_B)$. In particular, this implies that 
\begin{equation}\label{rho-tilde-rho} c_{\widetilde{C},F(B)}F(\widetilde{\rho}_B)=\rho_B,
\end{equation}
for any $B\in \Bc$. Using definition of $\Delta:C\to C\ot C$ given by diagramm \eqref{coalg-coend}, one can see that
\begin{equation}\label{D-d} c_{\widetilde{C}, C} F(\delta)=F(\Delta).
\end{equation}

If $h_1, h_2:C\ot F(B)\to C\ot C\ot F(B)$ are defined as
$$ h_1=\id_C\ot \rho_B, \quad h_2=\Delta\ot \id_{F(B)},$$
then, it is a standard result that, the kernel of the difference $h_1-h_2$ is $\rho_B$. Using   the above, we have the following:
\begin{prop}\label{exact-seq}  Assume that the module functor  $(F,c):\Bc\to \vc$ is faithful and  exact. The map $\widetilde{\rho}_B:B\to \widetilde{C} \triangleleft F(B)$ is the kernel of the difference of morphisms
\begin{equation*}
\xymatrix{
\widetilde{C} \triangleleft F(B ) \ar@<0.8ex>[r]^-{\,\, \widetilde{h}_1\,\,}\ar@<-0.8ex>[r]_-{ \widetilde{h}_2}&  \widetilde{C} \triangleleft (C\ot F(B ) ),
}
\end{equation*}
where $\widetilde{h}_1= \id_{\widetilde{C} } \triangleleft \rho_B$ and $\widetilde{h}_2= \delta\triangleleft  \id_{F(B)}.$
\end{prop}
\pf Note that, we are using that $F$ is right exact so that $F(\widetilde{C})=C$. See Lemma \ref{right-exact-on-coend}. Using \eqref{rho-tilde-rho}, \eqref{D-d} and the naturality of $c$ one can prove that
$$ c_{\widetilde{C}, C\ot F(B)} F(\widetilde{h}_1-\widetilde{h}_2)= (h_1-h_2) c_{\widetilde{C}, F(B)}.$$
This implies that $F (\widetilde{\rho}_B)=\ker F(\widetilde{h}_1-\widetilde{h}_2)=F(\ker \widetilde{h}_1-\widetilde{h}_2)$, since $F$ is also left exact. The result follows, since $F$ is faithful. 
\epf

\subsection{Reconstruction results}

The following Theorem is one of our main results and is a generalization of \cite[Thm. 2.2.8]{PS}.

\begin{teo}\label{recons-th} Let $\vc$ be a tensor category, and $\Bc$ be a  right  $\vc$-module category. 
Let $(F,c):\Bc\to \vc$ be an exact and faithful  module functor. Let $C=\oint^{B\in \Bc} F(B)\ot {}^*F(B)\in\vc$ be the relative coend coalgebra introduced in Proposition \ref{coalg-coend} and consider its category of comodules ${}^C\vc$. Then the functor $F$ factorizes into an equivalence of $\vc$-module categories $\widehat{F}: \Bc\stackrel{\sim}{\longrightarrow} {}^C\vc$  and the forgetful functor $\uc:{}^C\vc \to\vc$:
\begin{equation}
\xymatrix{
\Bc
\ar[rr]^-{\widehat{F}}
\ar[dr]_-{F}
&&
{}^C\vc
\ar[dl]^-{\uc}\ar@{}[d]|{}&\\
&\vc.&
}
\end{equation} 
\end{teo}
\pf Recall that we denote by
$$\pi_A:F(A)\ot {}^*F(A)\to \oint^{B\in \Bc} F(B)\ot {}^*F(B) $$
the dinatural transformations of the relative coend $C(F)$. The definition of the functor $\widehat{F}:\Bc\to {}^C\vc$ was given in Proposition \ref{factorization-F}. Since $F$ is faithful, $\widehat{F}$ is also faithful. We need to prove that $\widehat{F}$ is essentially surjective and full.

\medbreak

Let us prove first that $\widehat{F}$ is essentially surjective. Take $(M,\rho)\in {}^C\vc$. We endowed $C\ot M$ structure of left $C$-comodule via $\Delta$, then $\rho:M\to C\ot M$ is actually a morphism in ${}^C\vc$. We begin with the observation that, $(M,\rho)$ is the kernel (in ${}^C\vc$) of the difference of morphisms
\begin{equation*}
\xymatrix{
 C \ot M \ar@<0.8ex>[r]^-{\,\, \Delta\ot \id\,\,}\ar@<-0.8ex>[r]_-{ \id\ot \rho}&  C \ot C \ot M.
}
\end{equation*}
Recall from Proposition \ref{factorization-F} that, for any $A\in \Bc$, the object $F(A)$ has a left $H$-comodule structure given by $\rho_A=\psi^C_{F(A),F(A)}(\pi_A).$ Hence, for any $V\in \vc$, the object $F(A)\ot V$ has a left comodule structure given by $\rho_A\ot\id_V$. The next Claim will be crucial for the proof. To prove this claim, it is essential that the dinatural maps $\pi$ satisfy equation \eqref{dinat:coend:module:right}.

\begin{claim}\label{claim:cmodule} For any $A\in \Bc$, $V\in \vc,$ morphisms $c_{A,V}:F(A\triangleleft V)\to F(A)\ot V$ are $C$-comodule maps.
\end{claim}
\pf[Proof of Claim] We need to prove that
\begin{equation*}
(\psi^H_{F(A),F(A)}(\pi_A)\ot \id_V)c_{A,V}=(\id_H\ot c_{A,V}) \psi^H_{F(A\triangleleft V),F(A\triangleleft V)}(\pi_{A\triangleleft V})
\end{equation*}
This equation is equivalent to
\begin{equation}\label{c-comod-maps} 
(\pi_A\ot \id_{F(A)\ot V}) (\id \ot \coev_{F(A)}\ot \id_V) c_{A,V}= (\pi_{A\triangleleft V}\ot c_{A,V})(\id\ot \coev_{F(A\triangleleft V)} ).
\end{equation}
Since $\pi$ satisfies equation \eqref{dinat:coend:module:right}, it follows that
$$\pi_{A\triangleleft V}= \pi_A (\id_{F(A)}\ot {}^*(F(\id_A \triangleleft \coev_{V^*})) ) \big(\beta^{V^*}_{A\triangleleft V, A}  \big)^{-1}.$$
Using the description of the prebalancing $\beta$ given in \eqref{preb-s1}, we obtain that
\begin{align}\label{pi1}\begin{split} \pi_{A\triangleleft V}= &\pi_A (\id_{F(A)}\ot {}^*(c_{A\triangleleft V, V^*}  F(\id_A \triangleleft \coev_{V^*})) ) (\id\ot \phi^l_{F(A\triangleleft V),V^*})\\
&(c_{A,V}\ot \id_{{}^*F(A\triangleleft V)}).
\end{split} \end{align}
It follows from the naturality of $c$ that
$$c_{A, V\ot V^*} F(\id_A \triangleleft \coev_{V^*})=(\id_{F(A)}\ot \coev_{V^*}). $$
Using Equation \eqref{modfunctor11}, it follows that
$$c_{A\triangleleft V, V^*} F(\id_A \triangleleft \coev_{V^*})= (c^{-1}_{A,V}\ot \id_{V^*} )(\id_{F(A)}\ot \coev_{V^*}). $$
Whence
\begin{align}\label{pi2}\begin{split} 
\pi_{A\triangleleft V}&=\pi_A  (\id_{F(A)}\ot {}^*(\id_{F(A)}\ot \coev_{V^*})) (\id_{F(A)}\ot {}^*(c^{-1}_{A,V}\ot \id_{V^*} ) )\\
&(\id_{F(A)}\ot \phi^l_{F(A\triangleleft V),V^*})(c_{A,V}\ot \id_{{}^*F(A\triangleleft V)})\\
&=\pi_A  (\id_{F(A)}\ot {}^*(\id_{F(A)}\ot \coev_{V^*})) (\id_{F(A)}\ot \phi^l_{V^*,F(A)\ot V} )\\
&  ( c_{A,V}\ot {}^*(c^{-1}_{A,V}) ).
\end{split} 
\end{align}
The second equality follows from the naturality of $\phi^l$. It follows by a tedious, but straightforward, computation that
\begin{align}\label{coev-phi}\begin{split} \big({}^*(\id_{F(A)}\ot\coev_{V^*})\ot \id_{F(A)\ot V} \big)&\big( \phi^l_{V^*,F(A)\ot V}\ot\id_{F(A)\ot V}\big)\\&\big( \id_{V}\ot\coev_{F(A)\ot V}\big)= \coev_{F(A)}\ot\id_{V}
\end{split}
\end{align}
Let us prove now equation \eqref{c-comod-maps}.  The right hand side of \eqref{c-comod-maps} is equal to
\begin{align*}  &=(\pi_{A\triangleleft V} \ot\id_{F(A)\ot V} ) (\id\ot c_{A,V}) (\id_{F(A\triangleleft V)}\ot \coev_{F(A\triangleleft V)}) \\
&= \big(\pi_A \ot\id_{F(A)\ot V}  \big) \big(\id_{F(A)}\ot {}^*(\id_{F(A)}\ot \coev_{V^*})\ot\id_{F(A)\ot V}\big)\\
& \big(\id_{F(A)}\ot \phi^l_{V^*,F(A)\ot V}\ot \id_{F(A)\ot V} \big) \big(  c_{A,V} \ot {}^*(c^{-1}_{A,V})\ot \id_{F(A)\ot V} \big)\\
& (\id\ot c_{A,V}) (\id_{F(A\triangleleft V)} \ot\coev_{F(A\triangleleft V)})  \\
&=\big(\pi_A \ot\id_{F(A)\ot V}  \big) \big(\id_{F(A)}\ot {}^*(\id_{F(A)}\ot \coev_{V^*})\ot\id_{F(A)\ot V}\big)\\
& \big(\id_{F(A)}\ot \phi^l_{V^*,F(A)\ot V}\ot \id_{F(A)\ot V} \big) \big( c_{A,V}\ot \coev_{F(A)\ot V} \big)\\
&=(\pi_A\ot \id_{F(A)\ot V}) (\id_{F(A)} \ot \coev_{F(A)}\ot \id_V) c_{A,V}.
\end{align*}
The second equality follows from Equation \eqref{pi2}, the third equality follows from \eqref{duality-iso}, and the last equality follows from \eqref{coev-phi}. This finishes the proof of the Claim.
\epf

From this Claim follows that, natural isomorphisms $c$ endow the functor $\widehat{F}$ with structure of module functor.
Let $$\widetilde{C}=\widetilde{C}(F,c)=\oint^{B\in \Bc} B \triangleleft {}^*F(B)\in \Bc$$
be the object defined in Proposition \ref{existence-coend0}, together with themap
$$\delta:\widetilde{C}\to  \widetilde{C} \triangleleft C,$$
presented in Diagram \eqref{coproduct-tilde-c}.  Since $F(\widetilde{C})=C$, it also follows from Claim \ref{claim:cmodule} that,  there are $C$-comodule  isomorphisms
$$ F(\widetilde{C} \triangleleft  M )\simeq C \ot  M,$$
$$F(\widetilde{C} \triangleleft (C\ot M) )\simeq C\ot C \ot M.$$
Hence $(M,\rho)$ is the kernel of the difference of morphisms
\begin{equation*}
\xymatrix{
F(\widetilde{C} \triangleleft  M ) \ar@<0.8ex>[r]^-{\,\, F(h_1)\,\,}\ar@<-0.8ex>[r]_-{ F(h_2)}&  F(\widetilde{C} \triangleleft (C\ot  M) ),
}
\end{equation*}
where $h_1= \delta \triangleleft   \id_M, h_2=\id_{\widetilde{C}} \triangleleft   \rho$. Since $F$ is left exact $(M,\rho)\simeq F( \ker(h_1-h_2))$.  This proves that $\widehat{F}$ is essentially surjective. 

\medbreak

Let us prove now that $\widehat{F}$ is full. Take $A, B\in \Bc$ and $f:F(A)\to F(B)$ a $C$-comodule morphism. We have a commutative diagram
\begin{equation}\label{comodule}
\xymatrix{
 F(A)\ar[d]^{f}\ar[rr]^{\rho^A\quad}&& C\ot  F(A)  
 \ar[d]^{\id_H\ot f} \\
 F(B) \ar[rr]^{\rho^B }&&  C\ot  F(B).}
\end{equation}
 This implies that, we have a commutative diagram \begin{equation}\label{F-bef}
          \begin{split}
           \xymatrix{
         0 
         \ar[rr]  
         && 
         F(A)
         \ar[rr]^{ \rho_A} 
         \ar[d]^f 
         && 
         C \ot F(A )
		  \ar[rr]^{h^A_1-h^A_2}
		\ar[d]^{\id  \ot f } 
		&&
		C \ot (C\ot F(A ) )
		\ar[d]^{\id\ot f}
		\\
		 0
         \ar[rr] 
         && 
          F(B)
         \ar[rr]^{\rho_B} 
          && 
         C \ot F(B )
           \ar[rr]^{h^B_1-h^B_2}
			&&
			C \ot (C\ot F(B ) ).
\\
              }
           \end{split}
      \end{equation}
Where $h^A_1= \id_H\ot \rho_A, h^A_2=\Delta\ot \id_{F(A)}.$ Note that, by the universal property of the kernel, the map $f:F(A)\to F(B)$ is the unique morphism fitting in this diagram. Define
$$ \widetilde{C}=\widetilde{C}(F,c)=\oint^{B\in \Bc} B \triangleleft {}^*F(B)\in \Bc.$$
See Proposition \ref{existence-coend0} for the definition of this coend. Using Proposition \ref{exact-seq}, for any $B\in \Bc$, we have an exact sequence
$$0\to B\xrightarrow{\quad\widetilde{\rho}_B}
\widetilde{C} \triangleleft F(B ) \xrightarrow{\quad \widetilde{h}\quad }  \widetilde{C} \triangleleft (C\ot F(B ) ),
$$
where $\widetilde{h}_B=\id_{\widetilde{C} } \triangleleft \rho_B- \delta\triangleleft  \id_{F(B)}.$ For the definition of $\delta$ and $ \widetilde{\rho}_B$ see \eqref{coproduct-tilde-a}, \eqref{coproduct-tilde-c}. Since $f$ is a comodule morphism, by the universal property of the kernel, there exists a unique morphism $\eta$ that fits into the diagram

 \begin{equation}
          \begin{split}
           \xymatrix{
         0 
         \ar[rr]  
         && 
         A
         \ar[rr]^{ \widetilde{\rho}_A} 
         \ar[d]^\eta 
         && 
         \widetilde{C} \triangleleft F(A )
		  \ar[rr]^{\widetilde{h}_A}
		\ar[d]^{\id  \triangleleft f } 
		&&
		\widetilde{C} \triangleleft (C\ot F(A ) )
		\ar[d]^{\id \triangleleft (\id\ot f)}
		\\
		 0
         \ar[rr] 
         && 
          B
         \ar[rr]^{\widetilde{\rho}_B} 
          && 
          \widetilde{C} \triangleleft F(B )
           \ar[rr]^{\widetilde{h}_B}
			&&
			\widetilde{C} \triangleleft (C\ot F(B ) ).
\\
              }
           \end{split}
      \end{equation}
Applying $F$ to this diagram we obtain diagram \eqref{F-bef}. By the uniqueness of $f$, we get that $F(\eta)=f$, proving that $\widehat{F}$ is full. 
\epf

\subsection{Hopf algebras in $\vc$ constructed from a coend}

In the case $\Bc$ is a monoidal $\vc$-module category, and the functor $(F,c,\xi):\Bc\to \vc$ is a monoidal module functor, one can endow the coalgebra $C=C(F)$ with a Hopf algebra structure. In this Section, we shall explain the construction of the product and antipode on $C(F).$ The new hypothesis needed to construct such Hopf algebra is that $\vc$ posses a braiding.
\medbreak

To construct the product $m:C(F)\ot C(F)\to C(F)$ we shall use ideas from \cite{Majid}. In few words, we shall construct some natural module transformation $F\ot F\to C(F)\ot F\ot F$ and use Lemma \ref{a-omega} to find $m:C(F)\ot C(F)\to C(F)$.
\medbreak

\begin{lema}\label{def-d1} Let $\ca$ be a braided tensor category, $\Bc$ be a monoidal $\vc$-module category, and a monoidal module functor $(F,c,\xi):\Bc\to \vc$. Let us denote $ d_{B,V}$ the following natural isomorphism
$$d_{B,V}:{}^*F(B\triangleleft V)\to  {}^*V\ot{}^*F(B) , $$
$$ d_{B,V}= \big(\phi^l_{F(B), V} \big)^{-1} {}^*(c^{-1}_{B,V}).$$
Then, for any $B\in \Bc$, $V\in\vc$, we have
\begin{equation}\label{about-d} 
d_{B,{}^*V\ot V }=(\phi^l_{{}^*V,V}\ot \id)(\id_{{}^*V}\ot d_{B,{}^*V}) d_{B\triangleleft {}^*V, V },
\end{equation}
\end{lema}
\pf It follows by a straightforward calculation, using \eqref{modfunctor11}. \epf

Recall, from Proposition \ref{bimodule-funct}, the functors
$$H, \widetilde{H}:\Bc\boxtimes_\ca \Bc\to \ca,$$
determined by
$$H(A\boxtimes B)=F(A\ot B), \quad \widetilde{H}(A\boxtimes B)=F(A)\ot F(B),$$
for any $A, B\in \Bc$. Recall also that, we are denoting by $\pi_B:F(B)\ot {}^*F(B)\to C(F)$ the dinatural transformations, and the associated left $C(F)$-coaction
$\rho_B=\psi^C_{F(B),F(B)}(\pi_B).$

\medbreak

\begin{rmk}  In the case $\ca=\vect_\ku$, the following Proposition is trivial, since it only says that some natural transformation is additive. In the general case, where $\ca$ is arbitrary, it is far from obvious, and it is a  crucial step towards the  reconstruction of the product in $C(F)$. Its proof will highlight the importance of all required axioms.
\end{rmk}

For later use, $\pi$ has to satisfy Equation \eqref{dinat:coend:module:right}, using prebalancing \eqref{preb-s2}. This means that
\begin{equation}\label{pi-cf} \pi_{B\triangleleft V}=\pi_B \big(\id_{F(B)}\ot {}^*F(\id_B \triangleleft \coev_{V^*}) \big)\big(\id_{F(B)}\ot d^{-1}_{B\triangleleft V,V^* }\big)\big(c_{B, V}\ot \id_{{}^*F(B\triangleleft V)}\big).
\end{equation}

\begin{prop}\label{nat-bimod-pi} Let $\ca$ be a braided tensor category, $\Bc$ be a monoidal $\vc$-module category, and  $(F,c,\xi):\Bc\to \vc$ be a monoidal module functor. The natural transformation $\mu: \widetilde{H} \to C(F)\ot \widetilde{H}$ determined by the composition
$$\widetilde{H}(A,B)\xrightarrow{ \xi_{A,B}}  H(A, B)\xrightarrow{\rho_{A\ot B } } C(F)\ot  H(A, B)\xrightarrow{\id\ot  \xi^{-1}_{A,B}} C(F)\ot \widetilde{H}(A,B)$$
is a natural module transformation. That is $\mu\in \Nat_{\! m}( \widetilde{H}, C(F)\ot \widetilde{H})$.\end{prop}
\pf It follows from Proposition \ref{bimodule-funct} that $\xi$ is a natural module transformation. Once it has been established that the diagram
\begin{equation}\label{rho-mod-nat}
\xymatrix{
H(A, B \triangleleft V) \ar[d]^{e_{A,B,V}}\ar[rr]^{\rho_{A\ot B\triangleleft V}\quad}&& C(F)\ot  H(A, B\triangleleft V)
 \ar[d]^{\id_{C(F)}\ot e_{A,B,V} }\\
 H(A, B) \ot V\ar[rr]^{\rho_{A\ot B}\ot\id_V }&&  C(F)\ot  H(A, B)\ot V.}
\end{equation}
is commutative, the proof will follow. Recall, from the proof of Proposition \ref{bimodule-funct} that
$$e_{A,B,V}: F(A\ot (B\triangleleft V))\to F(A\ot B) \ot V,$$
$$e_{A,B,V}= c_{A\ot B, V} F(l^{-1}_{A\ot B, V}(\id_A\ot l_{B,V})) $$
We shall also keep the notation
$$\eta_{A,B,V}: A\ot (B \triangleleft V)\to (A\ot B) \triangleleft V,$$
$$\eta_{A,B,V}= l^{-1}_{A\ot B, V}(\id_A\ot l_{B,V}).$$  
Using the definition of $\rho_A$, one can see that diagram \eqref{rho-mod-nat} amounts to
\begin{align}\label{rho-mod-nat2} \begin{split}
&(\id_{C(F)}\ot e_{A,B,V})(\pi_{A\ot B \triangleleft V}\ot\id_{F(A\ot B \triangleleft V)})(\id_{F(A\ot B \triangleleft V)}\ot\coev_{F(A\ot B \triangleleft V)} )=\\
&=(\pi_{A\ot B}\ot\id_{_{F(A\ot B)}\ot V})(\id_{F(A\ot B)} \ot \coev_{F(A\ot B) }\ot\id_V)e_{A,B,V}
\end{split}
\end{align}
Using the dinaturality of $\pi$, one can see that the left hand side of \eqref{rho-mod-nat2} is equal to
\begin{align*}&=\big(\pi_{(A\ot B )\triangleleft V} (F(\eta_{B,V})\ot {}^*F(\eta^{-1}_{A,B,V}))\ot e_{A,B,V}\big)(\id_{F(A\ot B \triangleleft V)}\ot\coev_{F(A\ot B \triangleleft V)} )\\
&=\big(\pi_{(A\ot B )\triangleleft V} \ot c_{A\ot B,V}\big) \big(F(\eta_{A,B,V}) \ot  \coev_{F((A\ot B) \triangleleft V)}\big)\\
&=\big(\pi_{(A\ot B )\triangleleft V} \ot \id\big) \big(F(\eta_{A,B,V}) \ot \id\ot c_{A\ot B,V}\big) \big(\id\ot  \coev_{F((A\ot B) \triangleleft V)}\big)\\
&=\big(\pi_{(A\ot B )\triangleleft V} \ot \id\big) \big(F(\eta_{A,B,V}) \ot \id\ot c_{A\ot B,V}\big) \big(\id\ot  d^{-1}_{A\ot B,V}\ot c^{-1}_{A\ot B,V}\big)\\
&(\id_{F(A\ot B\triangleleft V)\ot {}^*V\ot \coev_{A\ot B}\ot \id_V }) ( \id_{F(A\ot B\triangleleft V)}\ot\coev_V)\\
&=\big(\pi_{(A\ot B )\triangleleft V} \ot \id\big) \big(F(\eta_{A,B,V}) \ot d^{-1}_{A\ot B,V}\ot \id_{F(A\ot B)\ot V}\big)\\
&(\id_{F(A\ot B\triangleleft V)\ot {}^*V}\ot \coev_{A\ot B}\ot \id_V )  ( \id_{F(A\ot B\triangleleft V)}\ot\coev_V)\\
&=(\pi_{A\ot B}\ot\id_{_{F(A\ot B)}\ot V})(\id_{F(A\ot B)} \ot \coev_{F(A\ot B) }\ot\id_V)e_{A,B,V}
\end{align*}
The second equality follows from the definition of $e_{A,B,V}$ and \eqref{duality-iso},  the fourth follows from \eqref{coevaluation-mod}. The last equation follows by using \eqref{pi-cf} and \eqref{about-d}.
\epf

\begin{teo}\label{hopf-rec} Assume $(\vc,\sigma)$ is a braided tensor category, $\Bc$ is a monoidal right $\vc$-module category, and $(F,c,\xi):\Bc\to \vc$ is an exact and faithful monoidal module functor with monoidal structure
$$\xi_{A,B}:F(A)\ot F(B)\to F(A\ot B). $$ 
The relative coend coalgebra $C(F)$ from Theorem \ref{recons-th} has an algebra structure, with unit $u=\pi_\uno$ 
and  product $m:C(F)\ot C(F)\to C(F)$ determined by
\begin{align}\label{product-C} \begin{split}&(m\ot \id)(\pi_A\ot  \pi_B\ot \id_{F(A)\ot F(B)})(\id\ot \sigma_{F(A), F(B)\ot {}^*F(B)}\ot \id_{F(B)})\\
&(\id_{F(A)}\ot \coev_{F(A)}\ot \id_{F(B)}\ot \coev_{F(B)})=(\id_{C(F)}\ot \xi^{-1}_{A,B})\rho_{A\ot B} \xi_{A,B},
\end{split}
\end{align}
for any $A,B\in \Bc$. The object $C(F)$ becomes a  bialgebra in $\vc$. 
Moreover, the equivalence of $\vc$-module categories  $\widehat{F}:\Bc\to {}^{C(F)}\vc$ is an equivalence of tensor categories.
\end{teo}
\pf  Recall from Lemma \ref{cc} that, if $f:{}^{C\ot C}\ca\to \ca$ is the forgetful functor, then $C(f)= C\ot C$ and the dinatural transformations of this coend are given by
$$\lambda_{(W,\rho)}:W\ot {}^*W\to  C\ot C,$$
$$\lambda_{(W,\rho)}= \bar{\psi}^{C\ot C}_{W,W}(\rho).$$
Since the functor $F:\Bc\to \vc$ factorizes as $F=\uc\circ \widehat{F}$, then the functor $\widetilde{H}=\ot\circ (F\boxtimes F): \Bc\boxtimes_\ca\Bc\to \ca$ factorizes as
\begin{equation*}
\vspace{.1cm}
\xymatrix{
\Bc\hphantom{_\mathcal{C}}\boxtimes_\mathcal{C} \Bc
\ar[rr]^{\widetilde{H}}
\ar[d]^{\widehat{F}\boxtimes \widehat{F}}
&&
\ca 
\\
{^C\mathcal{C}}\hphantom{_\mathcal{C}}\boxtimes_\mathcal{C} {^C\mathcal{C}}
\ar[d]^{\ot} 
&&
\\
{^{C\otimes C}\mathcal{C}
\ar[uurr]^{f}
\hphantom{^{C\otimes C}}}
&&
}
\end{equation*} 
Using Lemma \ref{equivalence-on-coend}, since $\widehat{F}$ is an equivalence of right $\ca$-module categories, it follows that 
$$C( \widetilde{H})=\oint^{X\in \Bc\boxtimes_\ca\Bc} \widetilde{H}(X)\ot {}^*\widetilde{H}(X)=C(F)\ot C(F).$$
Lemma \ref{equivalence-on-coend} also explains how to compute dinatural transformations of this coend. If $\pi_X: \widetilde{H}(X)\ot {}^*\widetilde{H}(X)\to C(F)\ot C(F)$, $X\in \Bc\boxtimes_\ca\Bc$, are the dinatural transformations associated to this coend, then
\begin{align*}\pi_{A\boxtimes B}&=\lambda_{(F(A)\ot F(B),\rho)}=\big(\id_{C(F)\ot C(F)}\ot \ev_D\big)\big(\id_C\ot\sigma_{F(A),C}\ot \id_{F(B)\ot {}^*D}\big)\\
& \big(\pi_A\ot \id_{F(A)}\ot \pi_B\ot\id\big) (\id_{F(A)}\ot \coev_{F(A)}\ot \id_{F(B)}\ot \coev_{F(B)}\ot \id_{{}^*D}),
\end{align*}
for any $A,B\in \Bc$. Here $D=F(A)\ot F(B)$. Here $\rho:F(A)\ot F(B)\to C(F)\ot F(A)\ot F(B) $ is the comodule structure of the tensor product according to formula \eqref{comodule-tensor-prod}. It follows from Lemma \ref{a-omega} that there is an isomorphism
\begin{equation*}
\begin{split} \omega: \Hom_\vc(C(F)\ot C(F),C(F))\to \Nat_{\! m}(\widetilde{H}, C(F)\ot \widetilde{H}),\\
\omega(g)_X= (g\ot \id_{\widetilde{H}(X)})\psi^{C(F)\ot C(F)}_{\widetilde{H}(X),\widetilde{H}(X)}(\pi_X).
\end{split}
\end{equation*} 
Since, by Proposition \ref{nat-bimod-pi}, $(\id_{C(F)}\ot \xi^{-1}_{A,B}) \rho_{A\ot B} \xi_{A,B}$ is a natural module transformation, that is
$$(\id_{C(F)}\ot \xi^{-1}_{A,B}) \rho_{A\ot B} \xi_{A,B}\in Nat_{\! m}(\widetilde{H}, C(F)\ot \widetilde{H}), $$
then there exists a morphism $m:C(F)\ot C(F)\to C(F)$ such that $$\omega(m)_{A\boxtimes B}=(\id_{C(F)}\ot \xi^{-1}_{A,B}) \rho_{A\ot B} \xi_{A,B}.$$ Using the rigidity axioms, the naturality of $\sigma$ and the formula for $\pi_{A\boxtimes B}$, this equation implies \eqref{product-C}. Using \eqref{product-C}, follows that
$$m\circ (m\ot \id)= m\circ (\id\ot m),$$
$$m \circ( u\ot \id )=\id=m\circ (\id\ot u).$$
\medbreak

It follows also from \eqref{product-C} that $\xi$ is a comodule morphism, giving the functor $\widehat{F}$ structure of monoidal functor. 
\epf 

\begin{defi}\label{definition:H(F)} For any monoidal module functor $F:\Bc\to \vc$ we shall denote by
$$H(F)=\oint^{B\in \Bc} F(B)\ot {}^*F(B)\in \vc $$
the bialgebra with product given by Theorem  \eqref{hopf-rec} and coproduct \eqref{coproduct-c}.
\end{defi}

\begin{rmk} In Theorem \ref{hopf-rec}, some hypothesis on the functor $F:\Bc\to \vc$ are superfluous. It follows from \cite[Corollaire 2.10]{De} that, if $F$ is right exact then it is exact and faithful. 
\end{rmk}

Let $(H,\Delta, m)$ be a bialgebra in $\vc$. In the next results, we shall be devoted to prove that the  bialgebra reconstructed in Theorem \ref{hopf-rec}, from the forgetful functor  $\uc:{}^H\vc\to\vc$, coincides with the original  bialgebra $H$. 

\medbreak

In Lemma \ref{cc} we already proved that the reconstructed coalgebra coincides with $H$. Moreover, if $(W, \rho_W)\in {}^H\vc$, then we have defined dinatural transformations
$$\pi_W:W\ot {}^*W\to  H,$$  $$\pi_W=\bar{\psi}^H_{W,W}(\rho_W)=(\id_H\ot \ev_W)(\rho_W\ot \id_{{}^*W}).$$ See \eqref{duals2} for the definition of isomorphisms $\bar{\psi} $. In order to see that, the reconstructed multiplication coincides with the product of $H$, we only need  to prove that the original product $m$  satisfies diagram \eqref{product-C}. This will be done in the next Proposition.

\begin{prop}  Let $(H,\Delta, m)$ be a  bialgebra in $\vc$. Using dinatural transformations $\pi_W:W\ot {}^*W\to  H,$ $\pi_W=\bar{\psi}^H_{W,W}(\rho_W),$ we have that equation
\begin{align*} m (\pi_V\ot \pi_W\ot \id_{V\ot W})&(\id_{V\ot {}^*V}\ot \sigma_{V, W\ot {}^*W}\ot \id_W)(\id\ot \coev_V\ot\id\ot \coev_W)= \\
&=\rho_{V\ot W}.
\end{align*}
holds for any pair $(V, \rho_V), (W,\rho_W)\in  {}^H\vc$.
\end{prop}
\pf It follows using the naturality of $\sigma$ and \eqref{comodule-tensor-prod}.
\epf

\subsection{The antipode of $H(F)$} In the next results we shall  construct an antipode $S:H(F)\to H(F)$, making the bialgebra $ H(F)$ a Hopf algebra in $\vc$.

\begin{lema}\label{antipode-lema}  Let $(\vc,\sigma)$ be a braided tensor category, and $\Bc$ be a monoidal  right  $\vc$-module category.
Let $(F,c):\Bc\to \vc$ be an exact and faithful  monoidal module functor.  For any $B\in \Bc$, $V\in \vc$ there are natural isomorphisms \begin{equation}\label{ts} t_{B,V}: {}^*(B \triangleleft V)\to {}^*B \triangleleft {}^*V\end{equation} such that the diagram
\begin{equation}\label{diagram-ts}
\xymatrix{
{}^*F(B \triangleleft V) \ar[d]^{{}^*(c^{-1}_{B,V})}\ar[rr]^{\qquad\; \qquad \; \qquad F(t_{B,V})\qquad \; \qquad\; \qquad}&& F({}^*B \triangleleft {}^*V)
 \ar[d]^{c_{{}^*B,{}^*V}} \\
 {}^*(F(B)\ot  V)\ar[rr]_{\qquad \qquad \; \qquad\sigma_{ {}^*V, {}^*F(B)} (\phi^l_{F(B),V})^{-1} \qquad \; \qquad\; \qquad}&&  F({}^*B) \ot {}^*V,}
\end{equation}
is commutative. 
\end{lema}
\pf  We shall freely use the fact that $F({}^*B)= {}^*F(B)$, for any $B\in \Bc$. Since we are under the same hypothesis as Theorem \ref{recons-th}, the functor $F:\Bc\to {}^H\vc$ is full. Since $F({}^*B \triangleleft {}^*V)$ is a left $H$-comodule, and the composition 
$$h=c^{-1}_{{}^*B,{}^*V} \sigma_{ {}^*V, F(B)} (\phi^l_{F(B),V})^{-1} {}^*(c^{-1}_{B,V}):{}^*F(B \triangleleft V)\to F({}^*B \triangleleft {}^*V),$$
 is an isomorphism, we can endow ${}^*F(B \triangleleft V)$ with some $H$-comodule structure such that, $h$ is a $H$-comodule map. Fullness of $F$ implies that, there exist some $t_{B,V}$ such that $F(t_{B,V})=h$.
\epf

We shall use the same notation as in previous sections. We denote by
$$\pi_A:F(A)\ot {}^*F(A)\to H(F)=\oint^{B\in \Bc} F(B)\ot {}^*F(B) $$
the dinatural transformations of the relative coend $H(F)$. Also, recall from Proposition \ref{factorization-F}, that for any $A\in \Bc$ we have that $(F(A),\rho_A)$ is a left $H(F)$-comodule, with structure given by
$$\rho_A=\psi^{H(F)}_{F(A),F(A)}(\pi_A)=(\pi_A\ot \id_{F(A)})(\id_{F(A)} \ot  \coev_{F(A)}).$$
Henceforth, for simplicity, we shall denote $H=H(F)$.
\begin{lema}\label{lemma-antip} Let $(\vc,\sigma)$ be a braided tensor category, and $\Bc$ be a monoidal  right  $\vc$-module category. For any $B\in \Bc$, $V\in \vc$, let us recall morphisms $ d_{B,V}$, defined in Lemma \ref{def-d1}, as
$$d_{B,V}:{}^*F(B\triangleleft V)\to  {}^*V\ot{}^*F(B) , $$
$$ d_{B,V}= \big(\phi^l_{F(B), V} \big)^{-1} {}^*(c^{-1}_{B,V}).$$
Then, for any $B\in \Bc$, $V\in\vc$, we have
\begin{equation}\label{about-rho1} (\sigma^{-1}_{ {}^*V, H}\ot \id_{ {}^* F(B)})(\id_H\ot d_{B,V})\rho_{ {}^*(B\triangleleft V )}=(\id_{{}^*V}\ot \rho_{{}^*B}) d_{B,V}.
\end{equation}
\end{lema}
\pf 
Using Lemma \ref{antipode-lema}, there are natural isomorphisms $t_{B,V}: {}^*(B \triangleleft V)\to {}^*B \triangleleft {}^*V$ such that
\begin{equation}\label{a-t} c_{{}^*B, {}^*V}F(t_{B,V})=\sigma_{{}^*V, F(B)} d_{B,V}.
\end{equation}
The naturality of $\phi^l$ and $c$ implies that
\begin{equation}\label{nat-phi9} 
F({}^*(\id_B\triangleleft \coev_V ))=({}^*(\coev_V)\ot \id_{{}^*F(B)}) d_{B,{}^*V\ot V },
\end{equation}
for any $B\in \Bc$, $V\in \vc$. Using the dinaturality of $\pi$ we obtain that
\begin{equation}\label{pi-and-t} 
\pi_{{}^*(B\triangleleft V)}= \pi_{{}^*B\triangleleft {}^*V}\big(F(t_{B,V})\ot {}^*F(t^{-1}_{B,V})  \big),
\end{equation}
 for any $B\in \Bc$, $V\in \vc$. Also, dinatural transformations $\pi$ satisfy \eqref{dinat:coend:module:right}, this implies that, for any $B\in \Bc$, $V\in \vc$, we have
 \begin{align}\label{pi-sp} \begin{split} \pi_{B \triangleleft {}^*V}= &\pi_B \big(\id_{F(B)}\ot F({}^*(\id_B\triangleleft \coev_V )) \big)    \big(\id_{F(B)}\ot d^{-1}_{B\triangleleft {}^*V,V} \big)\\
 &\big( c_{B,{}^*V} \ot \id_{{}^*F(B\triangleleft {}^*V)}\big) .
 \end{split}
 \end{align}
Next, when it is not absolutely necessary, as a space saving measure, we shall write the identities $\id$, without using subscripts. Using the definition or $\rho$ we obtain that 
$$(\sigma^{-1}_{ {}^*V, H}\ot \id_{ {}^* F(B)})(\id_H\ot d_{B,V})\rho_{ {}^*(B\triangleleft V )} d^{-1}_{B,V}$$ is equal to
\begin{align*} &=(\sigma^{-1}_{ {}^*V, H}\ot \id_{ {}^* F(B)})(\id_H\ot d_{B,V}) (\pi_{{}^*(B\triangleleft V)}\ot \id)(\id\ot \coev_{{}^*F(B\triangleleft V)})d^{-1}_{B,V}\\
&=(\sigma^{-1}_{ {}^*V, H}\ot \id_{ {}^* F(B)}) (\pi_{{}^*B\triangleleft {}^*V} \ot d_{B,V}) (F(t_{B,V})\ot {}^*F(t^{-1}_{B,V})\ot \id)\\
&(\id\ot \coev_{{}^*F(B\triangleleft V)})d^{-1}_{B,V}\\
&=(\sigma^{-1}_{ {}^*V, H}\ot \id_{ {}^* F(B)}) (\pi_{{}^*B\triangleleft {}^*V} \ot d_{B,V}) (F(t_{B,V})\ot \id_{{}^*F({}^*B \triangleleft {}^*V)} \ot  F(t^{-1}_{B,V}))\\
&(\id\ot \coev_{F({}^*B\triangleleft {}^*V)})d^{-1}_{B,V}\\
&=(\sigma^{-1}_{ {}^*V, H}\ot \id_{ {}^* F(B)})(\pi_{{}^*B}\ot \id) (\id\ot F({}^*(\id_{{}^*B}\triangleleft \coev_V ))\ot \id)(\id\ot d_{{}^*B \triangleleft {}^*V,V}\ot \id)\\
& (c_{{}^*B, {}^*V}\ot \id)\big(F(t_{B,V})d^{-1}_{B,V}\ot \id_{{}^*F({}^*B \triangleleft {}^*V)} \ot  d_{B,V} F(t^{-1}_{B,V})\big) (\id\ot d^{-1}_{{}^*B,{}^*V}\ot c^{-1}_{{}^*B, {}^*V} )\\
&(\id\ot \coev_{{}^*F(B)}\ot \id_{{}^*V})(\id_{{}^* F(B\triangleleft V)} \ot \coev_{{}^*V} )\\
&=(\sigma^{-1}_{ {}^*V, H}\ot \id_{ {}^* F(B)})(\pi_{{}^*B}\ot \id) (\id\ot ({}^*(\coev_V)\ot \id_{{}^*B}))d_{B,{}^*V\ot V }\ot \id)\\ &(\id\ot d_{{}^*B \triangleleft {}^*V,V}\ot \id)
 \big(c_{{}^*B, {}^*V} F(t_{B,V})d^{-1}_{B,V}\ot \id_{{}^*F({}^*B \triangleleft {}^*V)} \ot  d_{B,V} F(t^{-1}_{B,V})c^{-1}_{{}^*B, {}^*V}\big)\\ & (\id\ot d^{-1}_{{}^*B,{}^*V}\ot \id )
(\id\ot \coev_{{}^*F(B)}\ot \id_{{}^*V})(\id_{{}^* F(B\triangleleft V)} \ot \coev_{{}^*V} )\\
&= (\id_{{}^*V} \ot \pi_{{}^*B}\ot \id) (\sigma^{-1}_{ {}^*V, F({}^*B)\ot {}^*F({}^*B) }\ot \id)(\id \ot \ev_{{}^*V}\ot \id_{{}^*F({}^*B)\ot {}^*V\ot {}^*F(B)})\\
&(\sigma_{{}^*V, {}^*F(B)}\ot \id_{{}^{**}V\ot {}^*F({}^*B)} \ot \sigma_{{}^*V, {}^*F(B)}) (\id_{{}^*V\ot {}^*F(B)\ot  {}^{**}V}\ot \coev_{{}^*F(B)}\ot \id_{{}^*V})\\
&(\id\ot \coev_{{}^*V})\\
&=(\id_{{}^*V}\ot \rho_{{}^* B}).
\end{align*}
The second equality follows from \eqref{pi-and-t}, the third equality follows from \eqref{duality-iso}, the fourth one follows from \eqref{pi-sp}, the fifth equality follows from \eqref{nat-phi9}. The sixth equality follows by using \eqref{coev-of-dual} and \eqref{a-t}, and the last equality follows by using the rigidity axioms and \eqref{braiding1}.
\epf

Recall that, sometimes we are denoting $H=H(F)$.

\begin{teo}\label{antipode-natural}  Let $\vc$ be a tensor category, and $\Bc$ be a monoidal  right  $\vc$-module category.
Let $(F,c):\Bc\to \vc$ be an exact and faithful  monoidal module functor. There exists a map  $S:H(F)\to H(F)$, such that, it corresponds, under isomorphism $$\omega: \Hom_\vc(H(F),V)\to \Nat_{\! m}(F, V\ot F)$$ presented in \eqref{omega-isom}, to the natural module transformation $\alpha:F\to H(F)\ot F,$
\begin{align*}
\alpha_B=(\ev_{F(B)}&\ot \id_{H(F)\ot F(B)})(\id_{F(B)}\ot \sigma^{-1}_{{}^*F(B),H(F)}\ot \id_{F(B)})\\
& (\id_{F(B)}\ot \rho_{{}^*B}\ot \id_{F(B)})(\id_{F(B)}\ot \coev_{F(B)}).
\end{align*} 
That is $\omega(S)=\alpha$.
\end{teo}
\pf  We only need to prove that, indeed, $\alpha$ is a \textit{module} natural transformation, that is
$$(\id_H\ot c_{B,V}) \alpha_{B\triangleleft V}= (\alpha_B\ot\id_V) c_{B,V},$$
for any $B\in \Bc$, $V\in \vc$.  Equations \eqref{evaluation-mod} and \eqref{coevaluation-mod} implies that
\begin{equation}\label{ev-coev}\begin{split} \ev_{F(B\triangleleft V )} = \ev_{F(B)}\big(\id\ot \ev_V\ot \id \big)\big(c_{B,V}\ot d_{B,V}\big),\\
\coev_{F(B\triangleleft V )}= \big( d^{-1}_{B,V}\ot c^{-1}_{B,V}\big)\big(\id\ot \coev_{F(B)}\ot\id\big) \coev_V
\end{split}
\end{equation}

Using the definition of $\alpha$ we obtain that $(\id_H\ot c_{B,V}) \alpha_{B\triangleleft V} c^{-1}_{B,V}$ is equal to
\begin{align*} &= (\ev_{F(B\triangleleft V)}\ot \id_{H\ot F(B)\ot V)})(\id_{F(B\triangleleft V)}\ot \sigma^{-1}_{{}^*F(B\triangleleft V),H}\ot c_{B,V} )\\ &
 (\id_{F(B\triangleleft V)}\ot \rho_{{}^*(B\triangleleft V)}\ot \id_{F(B\triangleleft V)})(\id_{F(B\triangleleft V)}\ot \coev_{F(B\triangleleft V)})c^{-1}_{B,V}\\
 &=(\ev_{F(B)}\ot \id_{H\ot F(B)\ot V)})
 \big(\id\ot \ev_V\ot \id_{{}^*F(B) \ot H\ot F(B)\ot V} \big)\\
 &
 (\id_{F(B)\ot V}\ot \sigma^{-1}_{{}^*V\ot {}^*F(B),H}\ot \id_{F(B)\ot V} ) \big(c_{B,V}\ot  \id_{H}\ot d_{B,V} \ot c_{B,V}\big)\\
 &(\id_{F(B\triangleleft V)}\ot \rho_{{}^*(B\triangleleft V)}\ot \id_{F(B\triangleleft V)})\big( \id_{F(B\triangleleft V)}\ot d^{-1}_{B,V}\ot c^{-1}_{B,V}\big) \\ &\big(\id_{F(B\triangleleft V)\ot {}^*V}\ot \coev_{F(B)}\ot\id_V\big) (\id_{F(B\triangleleft V)}\ot \coev_V) c^{-1}_{B,V}\\
 &=(\ev_{F(B)}\ot \id_{H\ot F(B)\ot V)})
 \big(\id\ot \ev_V\ot \id_{{}^*F(B) \ot H\ot F(B)\ot V} \big)\\
 &(\id\ot \sigma^{-1}_{{}^*F(B),H}\ot \id)(\id\ot \sigma^{-1}_{{}^*V,H}\ot \id) \big(\id\ot (\id_{H}\ot d_{B,V} ) \rho_{{}^*(B\triangleleft V)} d^{-1}_{B,V} \ot \id \big)\\
 &\big(\id_{F(B\triangleleft V)\ot {}^*V}\ot \coev_{F(B)}\ot\id_V\big) (\id_{F(B\triangleleft V)}\ot \coev_V)\\
 &=\alpha_B\ot\id_V.
\end{align*}
The second equality follows from \eqref{ev-coev} and the naturality of $\sigma$, the third equality follows by \eqref{braiding1}. The last equality follows from \eqref{about-rho1} and rigidity axioms.
\epf

\begin{rmk} The above result is the most sensitive statement to prove in order to prove the existence of the antipode for $H(F)$. In the case $\vc=\vect_\ku$, the isomorphism
$$\omega: \Hom_\vc(H(F),V)\to \Nat(F, V\ot F)$$
lands in the space of \textit{all} natural transformations. So, in that case, there is nothing to prove, and the existence of $S$ is guaranteed by the fact that $\omega$ is an isomorphism.
\end{rmk}

\begin{cor}\label{antipode-reconstruct} The bialgebra $H(F)$, from Definition \ref{definition:H(F)}, is actually a Hopf algebra. The antipode $S:H(F)\to H(F)$ is determined as  the unique morphism such that
\begin{equation}\label{antipode}
\xymatrix{
 F(B)\ot {}^*F(B) \ar[d]^{\nu_B}\ar[rr]^{\qquad\pi_B}&& H(F)
 \ar[d]^{S} \\
  F(B)\ot {}^*F( B) \ot H\ar[rr]^{\qquad \ev_{F(B)} \ot\id_{H(F)} }&& H(F),}
\end{equation}
is commutative. Here 
$$\nu_B= (\id_{F(B)}\ot  \sigma^{-1}_{{}^*F(B),H(F)})(\id_{F(B)}\ot \pi_{{}^*B}\ot\id_{{}^*F(B)})(\id_{F(B)\ot {}^*F(B)}\ot \coev_{{}^*F(B)}).$$
\end{cor}
\pf Taking the map $S:H(F)\to H(F)$ defined in Theorem \ref{antipode-natural}, and using the definition of $\alpha$, one can prove that $S$ satisfies Diagram \eqref{antipode}. 

Axioms $$m\circ (\id_H\ot S)\circ \Delta=u\circ \varepsilon,$$
 $$m\circ ( S\ot\id_H)\circ \Delta=u\circ \varepsilon,$$
 follow from \eqref{antipode} by a lengthy, but straightforward, computation.
\epf

\begin{rmk} A version of the Hopf algebra $H(F)$ has been previously considered in the work of Lyubashenko \cite{Li1},  \cite{Li2}, when $F$ is some fiber functor and in \cite{Li3},  \cite{Li4}, when $F$ is the identity functor. See also the work of  P. Schauenburg \cite{PS}, where $\vc$ is the category of finite dimensional vector spaces. However, there are new ingredients in our definition.  We require that $\Bc$ is actually a monoidal $\vc$-module category, and the functor $F$ is a monoidal  module functor.  The use of the \textit{relative} coend is another new feature of our construction.
\end{rmk}

Let $(H,\Delta, m, S)$ be a  Hopf algebra in $\vc$. In the next Proposition we shall see that the reconstructed antipode given in Theorem \ref{antipode-natural}, out of the forgetful functor $\uc:{}^H\vc\to\vc$, coincides with the original  antipode $S:H\to H$.

\begin{prop}  Let $(H,\Delta, m, S)$ be a  Hopf algebra in $\vc$. Using the identification $H=H(\uc)$ proven in Lemma \ref{cc}, the  reconstructed antipode obtained in Corollary \ref{antipode-reconstruct} coincides with the antipode $S$ of $H$. 
\end{prop}
\pf Let us denote by $\widehat{S}:H\to H$ the map reconstructed in Proposition \ref{antipode-reconstruct}. For any $V\in \vc$, recall the isomorphism \eqref{omega-isom}
\begin{equation*}
\begin{split} \omega: \Hom_\vc(H(\uc),V)\to \Nat_{\! m}(\uc, V\ot \uc),\\
\omega(g)_{(W,\rho_W)}= (g\ot \id_{W})\psi^{H}_{W,W}(\pi_W),
\end{split}
\end{equation*}  
for any $(W,\rho_W)\in {}^H\vc$.  Corollary \ref{antipode-natural} implies that, if $\alpha:\uc\to H\ot \uc$ is the module natural  transformation defined as
\begin{align*}
\alpha_{(W,\rho_W)}=( \ev_{W}&\ot \id_{H\ot W})(\id_{W}\ot \sigma^{-1}_{{}^*W,H}\ot \id_{ W})\\
& (\id_{W}\ot \rho_{{}^*W} \ot\id_{W})(\id_W\ot \coev_W),
\end{align*}
for any $(W,\rho_W)\in {}^H\vc$, then 
$$\omega(\widehat{S})=\alpha. $$ 
 Using formula  \eqref{antipode-com} for $\rho_{{}^*W}$ we get that
$$\alpha_{(W,\rho_W)}= (S\ot \id_{W})\, \rho_W=\omega(S)_{(W,\rho_W)}, $$
for any $(W,\rho_W)\in {}^H\vc$. Whence $S=\widehat{S}$.
\epf

\section{Examples and Applications}\label{section_Example}

\subsection{Generic Example}

The next example was treated along the paper. At the end we proved that the reconstructed Hopf algebra from the forgetful functor $F:{^H}\vc\to \vc$ coincides with $H.$

\begin{exa}
Let $H$ be a Hopf algebra in a braided finite tensor category $\vc$ and consider the tensor category $\Bc={^H}\vc$. Then the forgetful functor admits an exact faithful monoidal functor $F:{^H}\vc\to \vc$ and  the counit $\varepsilon_H$ gives rise to a monoidal section ${^H}\vc\leftarrow \vc$ that turns ${^H}\vc$ into a $\vc$-module category. 
\end{exa}

\begin{exa}
Let $H\subset L$ be Hopf algebras in $\vect_\ku$ with $\iota:H\to L$ the inclusion and assume there is a left-inverse Hopf algebra map $\pi:H\rightarrow L$, called projection. Then we have the restriction functor $\iota^*\,:\Rep(H)\leftarrow \Rep(L)$ and in addition $\Bc=\Rep(L)$ becomes a modules category over $\vc=Rep(H)$ via $\pi^*:\,\Rep(H)\leftarrow \Rep(L)$. Assume that $H$ is quasitriangular. Hence our result shows that there exists a Hopf algebra $C\in \Rep(H)$ such that there is an equivalence of monoidal modular categories
$$\Rep(L)\cong {^C}\Rep(H).$$ 
We now discuss how this is related to the classical Radford Projection Theorem \cite{Rad85}. This result states that the existence of a Hopf algebra projection $\pi$ implies that $L$ is isomorphic to a Radford biproduct $$L\cong R\# H$$
where $R$ is a Hopf algebra in ${^H_H}\YD(\Vect_\ku)$. In our case we have assumed that $H$ is quasitriangular  and $\iota^*$ lifts to a braided functor to the center of $\Bc$. This gives rise to a choice of a functor $\Rep(H)\to {^H_H}\YD(\Vect_\ku)$ and shows $R$ to be in the image - more explicitly the $R$-matrix of $H$ determines the $H$-coaction from the $H$-action. Then this $R$ is precisely the dual of our Hopf algebra $C$.
\end{exa}

\subsection{Generic Consequences}

Our results has certain general implications, for example

\begin{exa}
If $(\vc, \sigma)$ is a braided tensor category, then the forgetful functor $f:Z(\vc)\to \vc$ has a section $G:\vc\to Z(\vc)$, $G(V)=(V, \sigma_V)$. Theorems \ref{recons-th} and  \ref{antipode-natural} imply that there exists a Hopf algebra $H\in \vc$ such that $Z(\vc)\simeq {}^H\vc$. For $\vc=\Rep(H)$ for a factorizable Hopf algebra $H$ this can be obtained from the defining equivalence  $Z(\vc)\simeq {}^H\vc\boxtimes {}^H\vc$, but in the non-nondegenerate case or the case without fibre functor, we are not aware of such a result. In general, we recover thereby the recent result \cite{LZ23}.
\end{exa}

\subsection{Lifting and cocycles deformations}
In the classification of (e.g. pointed) Hopf algebras $H$ the strategy in the Andruskiewisch-Schneider program \cite{AS10} is to consider the coradical $H_0$ (which is assumed to be a Hopf algebra, for example a group ring) and classify the possible Nichols algebras $R$, and then obtain $H$ as a lifting of $\mathrm{gr}(H)=R\#H_0$. It is an important question to determine whether this lifting is a $2$-cocycle twist, and the observation is that this holds in almost all cases \cite{AGI19}. \\

We view this problem in our setting: Since a main assumption is that  $H_0\subset H$ we have a tensor functor $\iota^*:\,\mathrm{Comod}(H_0)\to \mathrm{Comod}(H)$. The lifting is a $2$-cocycle twist precisely iff there is a tensor functor $\mathrm{Comod}(H_0)\leftarrow \mathrm{Comod}(H)$, if the $2$-cocycle is trivial then it comes from a Hopf algebra projection $\pi$, and then our relative coend is $R^*$. \\

While this view does not allow to decide the difficult  question when a lifting is a $2$-cocycle twist, it shows the natural categorical context of this question and it produces general statements, for example if $R$ as a coalgebra has a trivial lifting, then there is a monoidal section and hence the lifting is  tensor functor and thus the lifting comes from a $2$-cocycle twist. 

\subsection{The logarithmic Kazhdan Lusztig conjecture }

In conformal field theory, there is in good situations a modular tensor category of  representations of a vertex algebra $\mathcal{V}$, which reflects the analytic properties for example of solution spaces to certain differential equations (e.g. the braiding reflects the monodromy around the singularity $z=0$). One is often confronted with the very difficult question to determine the representation category of representations of a vertex subalgebra  $\mathcal{W}\subset \mathcal{V}$ if the representation theory of $\mathcal{V}$ is known, for example being a free-field realization. A brief introduction, references and an account for the statements below can be found in \cite{CLR}. \\

Categorically, one can understand as $\Rep(\mathcal{W})$ being a modular tensor category and $A=\mathcal{V}$ being a commutative  algebra in this category and $\vc=\Rep(\mathcal{V})$ is the category of local $A$-modules $\Rep^{loc}(A)$. Then by the results in \cite{CLR} Section $3$ we have under suitable conditions that $\Rep(\mathcal{W})$ is a relative Drinfeld center of $\Rep(A)$ with respect to the subcategory $\vc=\Rep^{loc}(A)$. In the cases relevant to logarithmic (i.e. nonsemisimple) conformal field theory, it is often the case that all \emph{simple} modules in $\Rep(A)$ are already in $\Rep^{loc}(A)$. This gives rise via the socle filtration to a monoidal section functor, see \cite{CLR} Section 4 $$\Rep^{loc}(A)\rightleftarrows \Rep^{loc}(A)$$
Then the results in our paper produce a Hopf algebra $R^*$, such that 
$$\Rep(A)={^{R^{*}}}\vc$$
and in this case the mentioned equivalence to the relative Drinfeld center means explicitly (see Section \ref{Section_YD})
$$\Rep(\mathcal{W})
\cong {^R_R}\YD(\vc)$$
For example the celebrated and notoriously difficult \emph{logarithmic Kazhdan Lusztig conjecture} considers the lattice vertex algebra $\mathcal{V}_\Lambda$ of the root lattice of a semisimple complex finite dimensional Lie algebra $\mathfrak{g}$ rescaled by an integer $p$, whose category of representations is the category of vector spaces graded by an abelian group 
$$\Rep(\mathcal{V}_\Lambda)=\Vect_{\Lambda^*/\Lambda}$$  Then it asserts that $\mathcal{V}_\Lambda$ contains as kernel of screening operators a certain vertex algebra $\mathcal{W}_p(\mathfrak{g})\subset \mathcal{V}_\Lambda$, whose category of representation is conjectured to be equivalent to the category of representations of the small quantum group $u_q(\mathfrak{g})$ for $q=e^{\pi i / p}$. In the first authors paper \cite{CLR} the approach is to view the small quantum group category as a category of Yetter-Drinfeld modules over the Nichols algebra $\mathcal{N}$
$$\Rep(u_q(\mathfrak{g}))\cong {^\mathcal{N}_\mathcal{N}}\YD(\Vect_{\Lambda^*/\Lambda})$$
and develop the categorical tools above to reduce the equivalence in question to proving that 
$$\Rep(A)\cong {_\mathcal{N}}\Vect_{\Lambda^*/\Lambda}$$
In the case $\mathfrak{sl}_2$ the abelian category could be determined and this gives a systematic proof of the conjecture in this case (initially we have used more complicated proof methods for the case $\mathfrak{sl}_2,p=2$ in \cite{CLR}, while \cite{GN} have used very different arguments for $\mathfrak{sl}_2,p$). A main motivation for writing the present paper is that now we have the clear statement that the category $\Rep(A)$ is for abstract reasons given as category of representations of a Hopf algebra $\mathcal{N}$ in $\vc$, and it now remains to determine that $\mathcal{N}$ is indeed the expected Nichols algebra, namely the Nichols algebra of screenings. 

\section{Questions}

\begin{question}
The result should be applicable if $\vc$ is merely locally finite tensor category an/or if $\mathcal{N}$ is infinite. The example we have in mind is maybe the quantum group $U_q(\mathfrak{g})$ at generic $q$.
\end{question}

\begin{question}
For given classes of semisimple modular tensor categories $\mathcal{C}$ (for example: the category of representations of an affine Lie algebra at positive integer level), can we classify semisimple Hopf algebras over $\mathcal{C}$?  
\end{question}

\begin{question}
Is there a general argument that the embedding of the coradical $H_0$ into a pointed Hopf algebra  $H$ admits a categorical section of comodule categories, if we assume in addition that the Yetter-Drinfeld module $H_1/H_0$ is semisimple over $H_0$?  
\end{question}

\end{document}